

\newif\ifabstract\abstractfalse		\newdimen\tempdim 
\newif\ifdraft\drafttrue		
\newif\ifslant\slantfalse		
\newif\ifthanks\thanksfalse

\newcount\definitionno 			\newbox\abstractbox
\newcount\chapterno 			\newbox\addressbox
\newcount\corollaryno 			\newbox\authorbox
\newcount\factno 			\newbox\mrbox
\newcount\lemmano 			\newbox\thankbox
\newcount \paragraphno 			\newbox\titlebox
\newcount\propositionno 
\newcount\sectionno \sectionno=-1
\newcount\subsectionno
\newcount\theoremno

\def\themonth{\ifcase \the\month \or Jan. \or Feb. \or March \or April
   \or May \or June \or July \or Aug. \or September \or Oct. \or Nov. \or
   Dec. \or Thermidor \fi}    
\def\thedate{\themonth \the \day, \the \year}

\def\abstract#1\par{\abstracttrue\def\theabstract{#1}}

\def\abstractbox{\ifabstract\vbox{\hsize 5 true in\leftskip 1 true in%
   \centerline{\sl Abstract}\smallskip%
   {\baselineskip 12 pt\noindent\sl\theabstract\par}}\vfill
   \ifthanks
   \centerline{\sl Acknowledgement}\smallskip
   \noindent\thethanks\par\message{Made abstract}
   \fi
   \else \message{No abstract supplied.}\fi}  

\def\thethanks{}
\def\thanks#1{\thankstrue\def \thethanks {#1}}
\def\title#1{\def\thetitle{#1}}

\def\authorbox{\vbox{\hbox{Saharon Shelah}  \medskip\hbox{\sl The
   Hebrew University} \hbox{\sl Rutgers University}}}
\def\titlebox{\vbox {\noindent\bigbf \thetitle}}

\def\firstpage{\pageno=0{\topskip= 2.5 true in \centerline
   {\vbox{\hskip 1 in\titlebox  \vskip.6true in \hskip 1 in\authorbox}}}
   \vskip 1 true in \abstractbox \vfill\eject}

\def\\#1\par{\medbreak\advance \paragraphno by 1 \medskip\noindent%
{\bf \sectionlabel \the\paragraphno \quad#1\par}\nobreak} 

\footline={\Sl \hfil
   \ifnum \pageno=0
      \ifdraft Draft of\/ {\bf\papername} \hfill \thedate\fi
      \else \ifnum\pageno<0 \rm \centerline{\folio}\fi
   \fi\hfil}  	
\def\draftheadline{\ifnum \pageno=0 \else\hfil \Sl Draft of {\em
   \papername}\quad -- \thedate\ -- 
   \quad\chaptername \quad\sectionname \hfil \rm \folio\fi}
   \headline={\ifdraft  \draftheadline \else \hfil 
	      \ifnum\pageno=0{}\else \rm \folio\fi 
           \fi}
\def\chaptername{
   \ifnum\chapterno>0Chapter \the \chapterno\quad\else
	\ifnum \pageno<0 Prefatory material\fi 
   \fi}
\def\proclaim#1{\ifvmode\ifx \xx#1\xx\else\bigskip\noindent \fi \fi
   {\bf #1}\par\nobreak\begingroup\Sl}
\def\endproclaim{\relax\par\endgroup\ifslant\errmessage{Font slanted
   after end of theorem}\fi\medskip} 
\def\proof#1:{\par\noindent{{\sl Proof}$\,$}#1:\par\penalty100} 
\def\references{\ifdraft \vfill \eject\else\bigbreak\fi
   \headline={\ifdraft \sl Draft of \rm \papername\hfil\quad \quad\hfil
   \sl References \qquad\else \hfil\fi\rm\folio}  \parindent 0 pt
   {\bf References\par\bigskip}\frenchspacing}
\newbox\refbox
\newdimen \refindent
\def\setrefwidth#1{\setbox\refbox=\hbox{[#1]xx}\refindent=\wd\refbox}
\setrefwidth{XXX}
\refindent=\wd\refbox
\def\journal#1**#2 (#3){{\sl #1} {\bf #2} (#3)}
\def\ref#1.{\bigskip\noindent \hangindent=\refindent\hangafter 1\hbox
   to \refindent{[#1]\hfil}\ignorespaces}  
\def\same{\vrule height.4pt depth 0pt width .6 in}


\def\JSL{\journal J. Symb. Logic**}
\def\LC#1,{{\sl in} {\bf Logic Colloquium #1},}
\def\LNM#1,{{\bf Lect. Notes Math. #1}, Springer, New York,}
\def\LNMb#1.{{\bf Lect. Notes Math. #1}, Springer, New York.}

\def\sectionname{\ifnum\sectionno>0 Section \the \sectionno\quad\else
\ifnum\paragraphno > 0 Par. \the\paragraphno\fi\fi}
\def\sectionlabel{\ifnum\sectionno>0 \the\sectionno.\fi}

\def\appendix{\ifdraft\vfill\eject\else \bigbreak\fi
\def\chaptername{Appendix}\def\sectionname{}\par\bigbreak\paragraphno=0\sectionno=0\def\sectionlabel{A}\noindent {\bf Appendix}\par\bigskip}
\def\introduction{\ifdraft\ifnum \pageno>0 \errmessage{Pageno should be set
negative during introduction}\fi\fi\bigskip\noindent{\bf Introduction}\medskip}

\def\section#1\par{\ifdraft\vfill\eject\else \bigbreak\fi\paragraphno=0\advance\sectionno by 1\noindent {\bf
\the\sectionno.\qquad  #1\par}\subsectionno=0\smallskip}



\font\bigbf=cmbx12

\input amssym.def \input amssym.tex 

\def\1{\rlap 1\kern.4pt1}
\def\<{\langle}
\def\>{\rangle}
\def\\{\qquad}
\def\={\equiv}
\def\and{\,\&\,}

\def\cite#1{{\rm[{\bf#1}]}}

\def\disjointunion{\rlap{$\union$}\kern1.8pt\raise.5em\hbox{$\cdot$}\kern2pt}

\def\implies{\Longrightarrow}
\def\includedin{\subseteq}
\def\intersect{\cap}
\def\iso{\simeq}

\def\normal{\triangleleft} 
\def\normalne{\raise1pt\hbox{\rlap{$\normal$}}_\ne}

\def\onto{\rlap{$\to$}\kern.05pt\to}

\def\midbar{\raise .4em\hbox{\vrule depth 0pt height .4pt width 4pt}\kern1pt}
\def\vtline{\kern 1 pt{\vrule depth 1pt height .8 em width .4pt}\kern.5pt}
\def\proves{\vtline\kern-.9pt\midbar}

\def\qed{\hfill{\vrule height 8pt depth .05pt width 5pt}}
\def\restrictedto{{\upharpoonright}}

\def\satisfies{\models}
\def\sdprod{{>}{\kern -2pt \triangleleft}}
\def\splitline#1//#2//#3{\vskip .1 in{\noindent\raise 10 pt\hbox{$#3$}\vbox{\hbox{#1}\hbox{#2}}\par}\vskip 4 pt}
\def\sqr#1#2{{\vbox{\hrule height .#2pt
	\hbox{\vrule width .#2pt height #1pt \kern#1pt
		\vrule width .#2pt}
		\hrule height .#2pt}}}
 
\def\th{^{\rm th\,}}
\def\to{\,\longrightarrow\,}

\def\with{\leftrightarrow}
\def \union{\cup}

\def\Disjointunion{\rlap{$\Union$}\kern2.7pt\raise.6em\hbox{$\cdot$}\kern2.5pt}

\def\Union{\bigcup}

\def\vline{\kern 1 pt{\vrule depth 1pt height .84 em width .4pt}\kern.5pt}


\def\fakebf#1{\rlap{$#1$}\kern.25pt\rlap{$#1$}\kern.25pt#1}
\def\bb#1#2{\rlap{\rm #1}\kern #2pt{\rm #1}}

\def\scriptbb#1{\rlap{$\rm \scriptstyle #1$}\kern.5pt\hbox{$\rm \scriptstyle#1$}}

\def\littleslant{\hbox{$\scriptscriptstyle/$}}
\def\slantline{\raise 1.5 pt \hbox{\rlap{\littleslant}\kern.1pt\raise.2pt\littleslant}}
\def\add#1#2#3{\rlap{\kern #3 pt #1}{\rm #2}}
	
\def\Ff{\Bbb F}
\def\Nn{\Bbb N}

\def\Rr{\Bbb R}
\def\Zz{\Bbb Z}

\def\AA{{\cal A}} 

\def\FF{{\cal F}} 

\def\MM{{\cal M}} 
\def\NN{{\cal N}} 
\def\PP{{\cal P}} 

\def\TT{{\cal T}}


\def\mod{{\rm \ mod\ }}

\def\tp{{\rm \ tp}}


\tolerance = 500

\def\oneandahalfspace{\baselineskip=18pt}
\def\singlespace{\baselineskip=\normalbaselineskip}


\def\Sl{\slanttrue\sl}

\def\em{\ifslant \rm \else \sl\fi}


\mathsurround=1pt
\parindent=25pt
\oneandahalfspace


\long\def\comment#1--{}

\comment 1:
 The control sequence: \\ generates a new, automatically numbered 
 paragraph, with the following material in boldface (until the next paragraph).
 For this to work properly, sections must begin with the \section command.
 The \section (and \chapter) commands also generate automatic numbering.

 While preparing drafts you should put ``\drafttrue'' at the
 beginning of the file.  The final version should start with
 ``\draftfalse''.
--

\comment 2:
   For the references we use two control sequences.
   The reference list is preceded by the command \references.
   Each reference has the form:
      \ref Label.   Author, title, ...
   and is followed by a blank line.
   This generates: [Label]   Author, title, ...
   with some good spacing.
--

\comment 3:
   Here is the file used to make paper 405.
--


\def\sss{\scriptscriptstyle}
\def\midbar{\raise .4em\hbox{\vrule depth 0pt height .4pt width 4pt}\kern1pt}
\def\up#1{\,{}^{#1}\kern -1.5pt}

\def\before{<\kern-5pt^{\sss|}_{\sss|}}
\def\cof{{\rm cof}\,}
\def\disjointunion{\rlap{$\union$}\kern1.8pt\raise.5em\hbox{$\cdot$}\kern2pt}
\def\diamond{\diamondsuit}
\def\dom{{\rm dom}\,}

\def\forces{\Vdash}
\def\td{\char'176}
\newbox\namebox\newdimen\tdwidth\newdimen\namewidth\setbox\namebox=\hbox{$\td$}
\tdwidth = \wd \namebox
\def\name#1{\setbox\namebox=\hbox{$#1$}\namewidth = \wd \namebox 
\advance \namewidth by - \tdwidth \divide \namewidth by 3
\kern \namewidth\lower 8 pt \hbox{\lower
\dp\namebox\rlap{\char'176}}\kern-\namewidth#1}      

\def\notforces{\rlap/\kern-2pt\forces}
\def\otp{{\rm otp\,}}
\def\phi{\varphi}

\def\uf{ultrafilter }

\font \tentt 	=cmtt10
\font \seventt  =cmtt8 at 7 pt
\font \fivett	=cmtt8 at 5 pt
\def\conditionfont{\textfont 1=\tentt\scriptfont 1=\seventt
   \scriptscriptfont 1=\fivett}
\def\condition#1{\hbox{$\conditionfont #1$}}
\def\pp{\condition  p}

\def\bold#1{{\bf #1}}
\def \pbar{\bold p}\def \qbar{\bold q} \def \rbar{\bold r}
\def \ubar{\bold u}
\def \xbar{\bold x}\def \ybar{\bold y}

\def\Hh{{\condition H}}
\def\Pp{{\condition P}}

\def\script#1{{\def\textstyle{\scriptstyle}#1}}



\def\bet{\beth}
\def\cf{\cof}

\def\epsilon{\varepsilon}
\def\eqdf{\mathrel{=:}}
\def\leend{\le_{\rm end}}

\def\App{{\it App }}

\def\Qset {{\condition Q}}
\def\Vof#1{V\left[{\script\Pp}\restrictedto #1\right]}

\newbox \commentbox \newcount \commentno
\def\comment#1{\global\advance \commentno by 1
   \setbox\commentbox \vbox{\hsize 2 true in \noindent\raggedright
   \tt\hangindent 20 pt  \hangafter 1 \the \commentno.\hskip 1 em 
   \strut #1 \strut}\leavevmode\vadjust{\smallskip \noindent \kern -20 pt \vbox{ \hrule
   \hbox {\vrule \kern 3 pt \box\commentbox \kern 3 pt\vrule}\hrule }}}

\def\offset#1{\hbox to \offsetsize {#1\hfil}}


\draftfalse

\ifdraft \singlespace \else \oneandahalfspace\fi
\ifdraft \pageno=-1\else \pageno = 1\fi

\def\papername{Paper 405}
\title {\vbox{\hbox{Vive la diff\'erence II.}
\hbox{The Ax-Kochen isomorphism theorem}}}

\abstract{We show in \S1 
that the Ax-Kochen isomorphism theorem [AK] requires the
continuum hypothesis. Most of the applications of this theorem are
insensitive to set theoretic considerations. (A probable exception is
the work of Moloney [Mo].) In \S2 we give an unrelated result
on cuts in models of Peano arithmetic which answers a question on the
ideal structure of countable ultraproducts of $\Zz$ posed in [LLS].
In \S1 we also answer a question of Keisler and Schmerl
regarding Scott complete ultrapowers of $\Rr$.}  

\thanks{The author thanks the Basic Research Fund of the Israeli  
Academy of Sciences, and the NSF
for partial support of this
research.\hfil
\break
\S1 of this paper owes its existence to
Annalisa Marcja's hospitality in Trento, July 1987; van den Dries'
curiosity about Kim's conjecture; and the willingness of Hrushovski and
Cherlin to look at \S3 of [326] through a glass darkly.
\S2 of this paper owes its existence to a question of G. Cherlin 
concerning [LLS].
This paper was prepared with the assistance of the group in Arithmetic
of Fields at the Institute for Advanced Studies, Hebrew University,
during the special year on Arithmetic of Fields, 1991-92. Publ. 405.}

\firstpage



\introduction

In a previous paper [Sh326] we gave two constructions of models of set
theory in which the following isomorphism principle fails in various
strong respects:
$$\vcenter{\hbox{If $\MM$, $\NN$ are countable elementarily equivalent
structures and $\FF$ is a nonprincipal}
\hbox{ultrafilter on $\omega$, then the ultrapowers $\MM^*$, $\NN^*$
of $\MM$, $\NN$ with  
respect to $\FF$ are isomorphic.}}\leqno (\hbox{Iso 1})$$
As is well known, this principle is a consequence of the continuum
hypothesis. 
Here we will 
give a related example in connection with the
well-known isomorphism theorem of Ax and Kochen. In its general
formulation, that result states that a fairly broad class of henselian
fields of characteristic zero satisfying a completeness (or
saturation) condition are classified up to isomorphism by the structure
of their residue fields and their value groups. The case that
interests us here is: $$\vcenter{\hbox{If $\FF$ is a nonprincipal
ultrafilter on $\omega$, then the ultraproducts}\hbox{$\prod_p
\Zz_p/\FF$ and $\prod \Ff_p[[t]]/\FF$ are isomorphic.}}
\leqno(\hbox{Iso 2})$$ Here $\Zz_p$ is the ring of  
$p$-adic integers and $\Ff_p$ is the finite field of order $p$.  
It makes no difference whether we work in the fraction fields of
these rings as fields, in the rings themselves as rings, or in the
rings as valued rings, as these structures are mutually interpretable
in one another.  In particular, the valuation is definable in the field
structure (for example, if the residual characteristic $p$ is greater than $2$
consider the property: ``$1+p x^2$ has a square root'').
We show that such an isomorphism cannot be obtained from the axioms of
set theory (ZFC).  As an application we may mention that certain
papers purporting to prove the contrary need not be refereed.

Of course, the Ax-Kochen isomorphism theorem is normally applied as a
step toward results which cannot be affected by set-theoretic
independence results. One exception is found in the work of Moloney
[Mo] which shows that the ring of convergent real-valued sequences on
a countable discrete set has exactly 10 residue domains modulo prime
ideals, assuming the continuum hypothesis. This result depends on the
general theorem of Ax and Kochen which lies behind the isomorphism
theorem for ultraproducts, and also on an explicit construction of a
new class of ultrafilters based on the continuum hypothesis. 
It is very much an open question to produce
a model of set theory in which Moloney's result no longer holds.

Our result can of course be stated more generally;
what we actually show here may be formulated as follows. 

\proclaim{Proposition A} It is consistent with the axioms of set theory
that there is an ultrafilter $\FF$ on $\omega$ such that for any two
sequences of discrete rank 1 valuation rings $(R^i_n)_{n=1,2,\ldots}$
($i=1,2$) having countable residue fields, 
any isomorphism $F:\prod_n R^1_n/\FF \to
\prod_n R^2_n/\FF$ is an ultraproduct of isomorphisms $F_n: R^1_n\to
R^2_n$ (for a set of $n$ contained in $\FF$). In particular most of
the pairs $R^1_n$, $R^2_n$ are isomorphic.
\endproclaim

In the case of the rings $\Ff_p[[t]]$ and $\Zz_p$,
we see that (Iso 2) fails.

From a model theoretic point of view this is not the right level of
generality for a problem of this type. There are three natural ways to
pose the problem:
$$\vcenter{\hbox{Characterize the pairs of countable models $\MM$,
$\NN$ such that for some ultrafilter $\FF$}
\hbox{in some forcing extension, $\prod \MM^\omega/\FF \not \iso \prod
\NN^\FF$;}} 
\leqno(1)$$ 
$$\vcenter{\hbox{Characterize the pairs of countable models $\MM$,
$\NN$ with no isomorphic ultrapowers}
\hbox{in some forcing extension;}}\leqno(2)$$ 
(there are two variants: the ultrapowers may be formed either
using one ultrafilter twice, or using any two ultrafilters). 
$$\vcenter{\hbox{Write $\MM\le \NN $ if in every forcing extension,
whenever $\FF$ is an ultrafilter}
\hbox{on $\omega$ such that $\NN^\omega/\FF$ is saturated, then
$\MM^\omega/\FF$ is also saturated. Characterize this relation.}}\leqno(3)$$ 
This is somewhat like the Keisler order [Ke, Sh-a or Sh-c Chapter VI] but does
not depend on the fact that the ultrafilter is regular. We can replace
$\aleph_0$ here by any cardinal $\kappa$ satisfying
$\kappa^{<\kappa}=\kappa$. 

However the set theoretic aspects of the Ax-Kochen theorem appear to
have attracted more interest than the two general problems posed here.
We believe that the methods used here are
appropriate also in the general case, but we have not attempted to go
beyond what is presented here. 

With the methods used here, we could try to show that for every $\MM$ with
countable universe (and language), if $\Pp_3$ is the partial order for
adding $\aleph_3$-Cohen reals then we can build a $\Pp_3$-name for a
non principal \uf $\FF$ on $\omega$, such that in $V^P$ $\MM^\omega/\FF$ 
resembles the models constructed in [Sh107]; we can choose the
relevant bigness properties in advance (cf. Definition 1.5, clause
(5.3)).  This would be helpful in 
connection with problems (1,2) above.

In \S2 of this paper we give a result on cuts in models
of Peano Arithmetic which has previously been overlooked. Applied to
$\omega_1$-saturated models, our result states that some cut does not
have countable cofinality from either side.  As we explain in \S2,
this answers a question on ideals in ultrapowers of $\Zz$ which was
raised in [LLS].  The result has nothing to do with the material in
\S1, beyond the bare fact that it also gives some information about
ultraproducts of rings over $\omega$.  

The model of set theory used for the consistency result in \S1 is
obtained by adding $\aleph_3$ cohen reals to a suitable ground model.  
There are two ways to get a ``suitable'' ground model. The first way
involves  taking any ground model which satisfies a portion of the GCH, 
and extending it by an appropriate
preliminary forcing, which generically adds the {\sl name} for an
ultrafilter which will appear after addition of the cohen reals.
The alternative approach is to start with an L-like ground model and
use instances of diamond (or related weaker principles) to prove that
a sufficiently generic name already exists in the ground model. 
That was the method used in \S3 of [Sh326],
which is based in turn on [ShHL162], which has still not appeared as
of this writing.  However the formalism of [ShHL162], though adequate 
for certain applications, turns out to be slightly too limited for our
present use.  More specifically, there are continuity assumptions
built into that formalism which are not valid here and cannot easily
be recovered. The difficulty, in a nutshell, is that a union of
ultrafilters in successively larger universes is not necessarily an
ultrafilter in the universe arising at the corresponding limit stages,
and it can be completed to one in various ways.  

We intend to include a more general version of
[ShHL162] in [Sh482].   
However as our present aim is satisfied by any model of set
theory with the stated property, we prefer to emphasize the first
approach here. So the family $\App$ defined below will be used as a
forcing notion for the most part. However we will also take note of
some matters relevant to the more refined argument based on a variant
of [ShHL162].  For those interested in such refinements, 
we summarize [ShHL162] in an appendix, as well as a version closer to
the form we intend to present in [Sh482].  
In addition the exposition in [Sh326, \S3] includes
a very explicit discussion of the way such a result may be used to
formalize arguments of the type given here, in a suitable ground model
(in the second sense).


\section Obstructing the Ax-Kochen isomorphism.

\\ Discussion

We will prove Proposition A as formulated in the introduction.  We
begin with a few words about our general point of view.  In practice
we do not deal directly with valuation rings, but with trees. If one
has a structure with a countable sequence of refining equivalence
relations $E_n$ (so that $E_{n+1}$ refines $E_n$) then the equivalence
classes carry a natural tree structure in which the successors of an
$E_n$-class are the $E_{n+1}$-classes contained in it. Each element of
the structure gives rise to a path in this tree, and if the
equivalence relations separate points then distinct elements give rise
to distinct paths. This is the situation in the valuation ring of of a
valued field with value group $\Zz$, where we have the basic
family of equivalence relations: $E_n(x,y)\iff v(x-y)\ge n$. 
(Or better: $E(x,y;z) =:$ ``$v(x-y) \ge v(z)$''.)
Of course an
isomorphism of structures would induce an isomorphism of trees, and
our approach is to limit the isomorphisms of such trees which are
available.

\\ The main result for trees.

We consider trees as structures equipped with a partial ordering and
the relation of lying at the same level of the tree. We will 
also consider expansions to much richer languages. 
We use the technique of [Sh326, \S3] to prove:

\proclaim{Proposition B} It is consistent with the axioms of set
theory that there is
a nonprincipal ultrafilter $\FF$ on $\omega$
such that for any two sequences of countable trees
$(T^i_n)_{n=1,2,\ldots}$ for $i=1,2$, 
with each tree $T^i_n$  countable with
$\omega$ levels, 
and with each node having at least two immediate successors,
if $\TT^i=\prod_n T^i_n/\FF$, 
then for any isomorphism $F:\TT^1\iso \TT^2$ there is an element $a\in
\TT^1$ such that the restriction of $F$ to the cone above $a$ is the
restriction of
an ultraproduct of maps $F_n:T^1_n\to T_n^2$.
\endproclaim

\\Proposition B implies Proposition A. 

Given an isomorphism $F$ between
ultraproducts $R^1$, $R^2$ 
modulo $\FF$ of discrete valuation rings $R^i_n$, we may consider 
the induced map $F_+$ on the 
tree structures $T^1$, $T^2$ associated with these rings, as 
indicated above. We then find 
by Proposition B that on a cone of $T^1$, $F_+$ agrees
with an ultraproduct
of maps $F_{+,n}$ between the trees $T^i_n$ associated with the $R^i_n$.
On this cone $F$ is definable from $F_+$, in the following sense:
$F(x) = y$ iff for all $n$, $F_+(a \mod \pi_1^n)\equiv b \mod \pi_2^n$,
where $\pi_i$ generates the maximal ideal of $R^i$ and we identify
$R^i/\pi_i^n$ with the $n$-th level of $T^i$. (This is expressed
rather loosely; in the notation we are using at the moment, one would
have to take $n$ as a nonstandard integer. After formalization in an
appropriate first order language it will look somewhat different.)
Furthermore $F$ is definable in $(R^1, R^2)$ 
from its restriction to this cone: the
cone corresponds to a principal ideal $(a)$ of $R^1$ and
$F(x)=F(ax)/F(a)$.
Summing up, then, there is a first order sentence valid in 
$(R^1,R^2; F_+)$ (with $F_+$ suitably interpreted as a parametrized
family of maps $R^1/\pi_1^n\to R^2/\pi_2^n$) stating that an
isomorphism $F:R^1\to R^2$ is definable in a particular way from $F_+$;
so the same must hold 
in most of the pairs $(R_{1,n},R_{2,n})$, that is, for a set of indices
$n$ which lies in $\FF$. 
In particular in such pairs we get an isomorphism of $R^1$
and  $R^2$. 

\\Context

We concern ourselves solely with Proposition B in the remainder of
this section. For notational convenience we fix
two sequences $(T^i_n)_{n<\omega}$ of
trees ($i=1$ or $2$) in advance, where each tree $T^i_n$ is countable with
$\omega$ levels, no 
maximal point, and no isolated branches.  The tree
$T^i_n$ is considered initially as a model with two relations: the tree order
and equality of level.
Although we fix the two sequences of trees,
we can equally well deal simultaneously with all possible pairs of
such sequences, at the cost of a little more notation.

As explained in the introduction, we work in a cohen generic
extension of a suitable ground model.
This ground model is assumed
to satisfy
$2^{\aleph_n} = \aleph_{n+1}$ for $n = 0, 1, 2$. 
If we use the partial
order \App defined below as a preliminary forcing, prior to the
addition of the cohen reals, then this is enough. If we wish to 
avoid any additional forcing then
we assume that the
ground model satisfies $\diamond_S$ for $S={\{\delta<\aleph_3:
\cf\delta=\aleph_2\}}$, 
and we work with \App directly in the ground model 
using the ideas of [ShHL162]. The second alternative requires more active
participation by the reader. 

Let $\Pp$ be cohen
forcing adding $\aleph_3$ 
cohen reals. 
An element $\pp$ of $\Pp$ is a
finite partial function from $\aleph_3\times\omega$ to $\omega$.  For
$\AA \subseteq\aleph_3$, and $\pp\in \Pp$, let $\pp\restrictedto \AA$
denote the restriction of $\pp$ to $\AA\times
\omega$ and $\Pp\restrictedto \AA=\{\pp\restrictedto \AA:\pp\in \Pp\}$. 
Let $\name x_\beta$ be the $\beta\th$ 
cohen real. The partial order \App is defined below.

We will deal with a number of expansions of the basic language of pairs of
trees. For a forcing notion $\Qset$ and $G$ $\Qset$-generic over $V$, we
write $^G(T^1_n,T_n^2)$ for the expanded 
structure in which for every $k$, every sequence
$(r_n)_{n< \omega}$ of $k$-place relations $r_n$ on $(T_n^1,T_n^2)$ is
represented by a $k$-place relation symbol $R$ (i.e.,
$R_{(r_n:n<\omega)}$); that is, $R$ is interpreted in 
$(T^1_n,T_n^2)$ by the relation $r_n$. This definition takes place in
$V[G]$. In $V$ we will have names for these relations and relation
symbols. We 
write $^\Qset(T^1_n,T_n^2)$ for the corresponding collection of names.  In
practice $\Qset$ will be $\Pp\restrictedto \AA$ for some
$\AA\includedin \omega_3$ and in this case we write
$^\AA(T_n^1,T_n^2)$. 

Typically we will have certain
subsets of each $T^i_n$ singled out, and we will want to study the
ultraproduct of these sets, so we will make use of the predicate
whose interpretation in each $T^i_n$ is the desired set. We would prefer
to deal with $^\Pp(T^1_n,T_n^2)$, but this is rather large, and so
we have to pay some attention to matters of timing. 

\\ Definition

As in [Sh326], we set up a class \App of approximations to the name of
an ultrafilter in the generic extension $V[\Pp]$.
In [Sh326] we emphasized the use of the 
general method of [ShHL162] to construct the name $\name \FF$ of a suitable 
ultrafilter in the ground model.  Here we emphasize the alternative
and easier approach, forcing with \App.
However we include a summary of the formalism of [ShHL162], and a
related formalism, in an
appendix at the end.

The elements of $\App$ are triples $q=(\AA, \name\FF, \epsilon)$ such that:
$$\hbox{$\AA$ is a subset of $\aleph_3$ of cardinality
   $\aleph_1$;}\leqno(1)$$ 
$$\hbox{$\name \FF$ is a $\Pp\restrictedto \AA$-name of a nonprincipal
   ultrafilter on $\omega$, called $\name \FF\restrictedto \AA$;}\leqno(2)$$ 
$$\hbox{$\epsilon=(\epsilon_\alpha:\alpha\in \AA)$, with each $\epsilon_\alpha\in
   \{0,1\}$, 
   and $\epsilon_\alpha=0$ whenever $\cf \alpha<\aleph_2$;}
   \leqno(3)$$

$$\hbox{For $\beta\in \AA$ we have: [$\name \FF\cap \{\name a:\name a$ a
   $\Pp\restrictedto (\AA\cap\beta)$-name of a subset of $\omega\}$] is a  
   $\Pp\restrictedto (\AA\cap\beta)$-name;}\leqno(4)$$
$$\hbox{If $\cf \beta=\aleph_2$, $\beta\in\AA,
   \epsilon_\beta=1$ then $\Pp\restrictedto \AA$ forces the 
   following:}\leqno(5)$$
\item{(5.1)} 
   $\name x_\beta/\name \FF$ is an element of 
   $(\prod_{n<\omega } T_n^{1}/\name \FF\restrictedto \AA)
   ^{\Vof \AA}$
   whose level is above all levels of elements of the form 
   $\name x/\name\FF$ for 
   $\name x$ a $\Pp\restrictedto(\AA\cap\beta)$-name;
 \item{(5.2)} $\name x_\beta$ induces a branch $\name B$ on
    $(\prod_n T^{1}_n)^{\Vof{(\AA\cap\beta)}}
    /[\name \FF\restrictedto (\AA\cap\beta)]$
    which has elements in every level of that tree  (such a branch
will be called {\em full}) and which is a
    $\Pp\restrictedto(\AA\cap\beta)$-
    name (and not just forced to be equal to one);
\item{(5.3)} The branch $\name B$ intersects every dense subset of
    \quad $(\prod_n^{\AA \intersect \beta}
    T^{1}_n)^{\Vof{(\AA\cap\beta)}} 
    /[\name \FF\restrictedto (\AA\cap\beta)]$
   which 
is definable in 
    $(\prod_n {}^{\AA \intersect \beta}(T^1_n,T_n^2)/[\name \FF\restrictedto
   (\AA\cap\beta)])^{\Vof{(\AA\cap\beta)}}$. 

Note in (5.3) that the dense subset under consideration will have a 
$\Pp\restrictedto (\AA\cap\beta)$-name, and also that 
by \L o\'s' theorem a dense subset of the type described 
extends canonically to a dense subset in any larger model.
The notion of ``bigness'' alluded to in the introduction is given by
(5.3).

We write $q_1\le q_2$ if $q_2$ 
extends $q_1$ in the natural sense. 
We say that $q_2\in \App$ is an {\em end} extension of $q_1$, and we
write $q_1\leend q_2$, if $q_1\le q_2$ and
$\AA^{q_2}\setminus\AA^{q_1}$ follows $\AA^{q_1}$.
Here we have used the notation: $q=(\AA^q, \name\FF^q,\epsilon^q)$. 

\\ Remark

The following comments bear on the version based on the method of [ShHL162].
In this setting, rather than examining each $\name x_\beta$ 
separately, we would
really group them into short blocks
$X_\beta=(\name x_{\beta+\zeta}:\zeta<\aleph_2)$, for $\beta$ divisible by
$\aleph_2$.  Then
our assumptions on the ground
model $V$ allow us to use the method of 
[ShHL162] to construct the name $\name \FF$
in $V$.  One of the ways $\diamond_S$ would be used is to ``predict''
certain elements $\pp_\delta\in \Pp\restrictedto \delta$ and certain
$\Pp\restrictedto \delta$-names of functions $\name F_\delta$ which
amount to guesses as to the 
restriction to a part of $\prod_n T^1_n$ of (the name of) a function
representing some isomorphism $\name F$ modulo $\name \FF$.
As we indicated at the outset, we intend to elaborate on 
these remarks elsewhere.

\\  Lemma

\proclaim{}
If $(q_\zeta)_{\zeta<\xi}$ 
is an 
increasing sequence of at most $\aleph_1$ members of $\App$ such that
$q_{\zeta_1}\leend q_{\zeta_2}$  for $\zeta_1<\zeta_2$,
then we can find $q\in \App$ such that $\AA^q=\bigcup_\zeta
\AA^{q_\zeta}$ and  $q_\zeta\leend q$ for $\zeta<\xi$.
\endproclaim

\proof:

We may suppose $\xi>0$ is a limit ordinal.
If $\cf(\xi)>\aleph_0$ then $\bigcup_{\zeta<\xi}q_\zeta$ will do,
while if $\cf(\xi)=\aleph_0$ then we just have to extend
$\bigcup_\zeta \name \FF^{q_\zeta}$ to a
$\Pp\restrictedto(\bigcup_\zeta\AA^{q_\zeta})$-name
of an ultrafilter on $\omega$, which is no problem.
(cf. [Sh326, 3.10]). 
\qed

\\Lemma 

\proclaim{}
Suppose $\epsilon = 1$, 
$q\in\App$,
$\gamma>\sup \AA^q$, 
and $\name B$ is a $\Pp\restrictedto \AA^q$-name
of a branch of $(\prod_n T^\epsilon_n/\name \FF^q)^{V[\Pp\restrictedto \AA^q]}$.
Then:
\item{1.} We can find an $r\in\App$ with  $\AA^r=
\AA^q\cup\{\gamma\}$, and a $(\Pp\restrictedto \AA^r)$-name $\name x$ of a member
of $\prod_n  T^\epsilon_n/\name \FF^r$ which is above $\name B$.
\item{2.} We can find an $r\in\App$ with $q\leend r$ and $\AA^r=
\AA^q\cup[\gamma,\gamma+\omega_1)$, and a $(\Pp\restrictedto \AA^r)$-name
$\name B'$ of a full branch extending $\name B$, which intersects
every definable dense subset of  
$(\prod_n {}^{\AA^r}T^\epsilon_n)^{V[\Pp\restrictedto \AA^r]}/\name \FF^r$.    

\item{3.} In (2) we can ask in addition that any particular type $p$
over $\prod {}^{\AA^q}(T_n^1,T_n^2)/\name \FF^q$ (in $V[\Pp\restrictedto \AA^q]$)
be realized in $(\prod_n {}^{\AA^r}T^\epsilon_n)^{V[\Pp\restrictedto
\AA^r]}/\name \FF^r$.      
\endproclaim

\proof:
1.  Make $\name x_\gamma$ realize the required type, and let
$\epsilon_\gamma=0$.

2. We define $r_\zeta=r\restrictedto (\AA^q\cup[\gamma, \gamma+\zeta))$
by induction on $\zeta\le \omega_1$. For limit $\zeta$ use 1.7 and
for successor $\zeta$ use part (1). One also takes care, via 
appropriate bookkeeping, that $\name B'$ should intersect every dense
definable subset of ($\prod_n{}^{\AA^r}T^\epsilon_n/\name \FF^r)^
{\Vof{\AA^r}}$ by arranging for each
such set to be met in some specific 
$(\prod_n{}^{\AA^{r_\zeta}}T^\epsilon_n/\name \FF^{r_\zeta})
^{V[\Pp\restrictedto \AA^{r_\zeta}]}$ with $\zeta<\aleph_1$.

3. We can take $\alpha \in [\gamma, \gamma+\omega_1)$ with $\cof \alpha 
\ne \aleph_2$ and use $x_\alpha$ to
realize the type.
\qed

\\ Lemma

\proclaim{}
Suppose $q_0, q_1,q_2\in \App$,
$q_0=q_2\restrictedto \beta$,
$q_0\le q_1$,
$\AA^{q_1}\subseteq\beta.$

\item{1.} 
If $\AA^{q_2}\setminus\AA^{q_0}=\{\beta\}$ and 
$\epsilon^{q_2}_\beta=0$, then
there is $q_3\in \App$, $q_3\ge q_1, q_2$
with $\AA^{q_3}=\AA^{q_1}\cup\AA^{q_2}$. 

\item{2.} 
Suppose $\AA^{q_2}\setminus \AA^{q_0}=\{\beta\}$,
$\cf \beta=\aleph_2$, $\epsilon^{q_2}_\beta\ne 0$, and in particular
$\sup \AA^{q_1}<\beta$.
Assume that
$\name B_1$ is a $\Pp\restrictedto \AA^{q_1}$-name of a full branch of
$(\prod T_n^{\epsilon^{q_2}_\beta}/\name \FF^{q_1})^{V[\Pp
\restrictedto \AA^{q_1}]}$ intersecting  
every dense subset of this tree which is definable in 
$(\prod_n {}^{\AA^{q_1}}(T_n^1,T_n^2)/\name \FF^{q_1})^
{\Vof{\AA^{q_1}}}$,    
such that $\name B_1$
contains the branch $\name B_0$ which $\name x_\beta$ induces 
according to $q_2$.
Then there is $q_3\ge q_1,q_2$ with $\AA^{q_3} =
\AA^{q_1}\cup\{\beta\}$, 
such that
according to $q_3$, $\name x_\beta$ induces 
$\name B_1$ on $(\prod T^{\epsilon^{q_2}_\beta}_n/\name
\FF\restrictedto \AA^{q_1})^{\Vof{\AA^{q_1}}}$.

\item{3.} If $\AA^{q_2}\setminus \AA^{q_0}=\{\beta\}$,
$\cf \beta=\aleph_2$, $\epsilon_\beta^{q_2}=1$, and 
$\sup\AA^{q_1}<\gamma<\beta$ with $\cf\gamma\neq \aleph_2$,
then there is $q_3\in \App$ with $q_1\le q_3$, $q_2\le q_3$,
$\AA^{q_3}=\AA^{q_1}\cup\AA^{q_2}\cup[\gamma, \gamma+\omega_1)$.
\item{4.} 
There are $q_3\in \App$, $q_1, q_2\le q_3$, 
so that $\AA^{q_3}\setminus\AA^{q_1}\cup\AA^{q_2}$ has the
form $\Union \{[\gamma_\zeta, \gamma_\zeta+\omega_1): \zeta \in
\AA^{q_2}\setminus\AA^{q_0}$, $\cf\zeta=\aleph_2\}$
where $\gamma_\zeta$ is arbitrary subject to 
$\sup(\AA^{q_2}\restrictedto \zeta)<\gamma_\zeta<\zeta$.
\item{5.} Assume $\delta_1 < \aleph_2$, $\beta < \aleph_3$, 
that $(p_i)_{i< \delta}$ is an
increasing sequence from \App, and that $q\in \App\restrictedto \beta$ 
satisfies:
$$
\hbox {For $i<\delta_1$: $p_i\restrictedto \beta\le q$.}
$$
Then there is an 
$r\in \App$    with 
$q \leend r$ and $p_i\le r$    for all    $i< \delta_1$. 
\item{6.} Assume $\delta_1, \delta_2 < \aleph_2$,
$(\beta_j)_{j<\delta_2}$
is an increasing sequence with all $\beta_j  < \aleph_3$, that
$(p_i)_{i< \delta_1}$ is an 
increasing sequence from \App, and that $q_j\in \App\restrictedto \beta_j$ 
satisfy:
$$
\hbox {For $i < \delta_1$, $j <\delta_2$: $p_i\restrictedto \beta_j\le q_j$;}
\qquad
\hbox {For $j < j' <\delta_2$: $q_j\leend q_{j'}$.}
$$
Then there is an 
$r\in \App$    with $p_i\le r$     and 
$q_j \leend r$ for all    $i< \delta_1$ and $j < \delta_2$. 
\endproclaim
\proof:
1. The proof is easy and is essentially contained in the proofs following.
(One verifies that $\name \FF^{q_1}\union \name \FF^{q_2}$ generates a proper
filter in $V[\Pp\restrictedto (\AA^{q_1}\union \AA^{q_2})]$.) 

2. Let $\AA_i=\AA^{q_i}$ and let $\name \FF_i=\name \FF^{q_i}$ for
$i=1,2$, and $\AA_3=\AA_1\union\AA_2=\AA_1\union\{\beta\}$.
The only nonobvious part is to show that in 
$\Vof{\AA_3}$
there is an ultrafilter extending $\name \FF_1 \union\name \FF_2$
which contains the sets: 
$$\{n:T_n^1\satisfies \name x(n)\le \name x_\beta(n)\}\;\hbox{ for
$\name x\in \name B_1,\;\name x$ 
a $\Pp\restrictedto \AA_1$-name.}$$ 
If this fails, then there is some $\pp\in \Pp\restrictedto \AA_3$,
a $\Pp\restrictedto \AA_1$-name $\name a$ 
of a member of $\name \FF_1$,
a $\Pp\restrictedto \AA_2$-name $\name b$ 
of a member of $\name \FF_2$,
and some $\name x\in \name B_1$
such that $\pp \forces$``
$\name a\cap \name b \cap\name c =\emptyset$''
where 
$\name c=\{n:\name x(n) \le \name x _\beta(n)\}$.
Let $\pp_i=\pp\restrictedto \AA_i$ for 
$i=0,1,2$, and let $\Hh^0\subseteq\Pp\restrictedto \AA_0$ be generic over $V$,
with $\pp_0\in \Hh^0$.

Let: 
$$\name A^1_n[\Hh^0]=
\{y\in T^1_n:
\hbox{ For some $\pp'_1, \pp_1\le \pp'_1\in \Pp\restrictedto \AA_1$, 
$\pp'_1\restrictedto \AA_0\in \Hh^0$ and 
$\pp'_1\forces$ 
``$\name x(n)\le y$, and $n\in \name a $''}\}.
$$
Then $\name A^1_n$ is a $\Pp\restrictedto \AA_0$-name.
Let $\name A^{1}=(\prod_n \name A^1_n/\name \FF\restrictedto
\AA_0)^{V[\Pp\restrictedto \AA_0]}$. 
Now $\name A^{1}$ is not necessarily dense in 
$(\prod_n T^1_n/\FF\restrictedto \AA_0)^{\Vof{\AA_0}}$, 
but the set 
$$
\name A^*\eqdf
\{\name y\in (\prod_n {}^*T^1_n 
/\name \FF^{q_0})^{V[\Pp\restrictedto \AA_0]}: 
\name y\in A^{1}, \hbox{ or $\name y$ is incompatible in the tree with 
all $\name y'\in A^{1}$}\}
$$ 
{\em is} dense, and it is definable, hence not disjoint from $\name B_0 $.
Fix $\name y\in \name A^*\intersect \name B_0$. 
As $\name x\in \name B_1$, $\name x$ and $\name y$ cannot be
forced to be incompatible, and thus $\name y\in \name A^{1}$. 

The following sets are in $\name \FF^{V[\Hh^0]}$:
$$\name A=\{n:\hbox{ for some }\pp'_1,
\pp_1\le \pp'_1\in \Pp\restrictedto \AA_1, \pp'_1\restrictedto \AA_0\in \Hh^0
\hbox{ and }\pp'_1\forces\hbox{``$\name x(n)\le\name y(n)$, and $n\in
\name a$''}\}.$$
$$\name B=\{n:\hbox{ for some }\pp'_2,
\pp_2\le \pp'_2\in \Pp\restrictedto \AA_2, \pp'_2\restrictedto \AA_0\in \Hh^0
\hbox{ and }
\pp'_2\forces\hbox{``$\name y(n)\le \name x_\beta(n)$, and $n\in \name
b$''}\}.$$ 
For example, $\name A$ is a subset of $\omega$ in $V[\Hh^0]$ which is
in $\name \FF^{q_1}$. As the complement of $\name A$ cannot be 
in $\name \FF^{q_0}$, $\name A$ must be.

Now for any $n\in \name A\intersect \name B$ we can force $n\in \name
a\intersect\name b\intersect\name c$ by amalgamating the corresponding
conditions $\pp'_1,\pp'_2$.

3. Let $\name B_0$ be the $\Pp\restrictedto \AA^{q_0}$-name of the
branch which $\name x_\beta$ induces.  By 1.8 (2) there is $q^*_1$,
$\AA^{q^*_1}=\AA^{q_1}\cup[\gamma, \gamma+\omega_1)$, 
$q_1\le q^*_1\in \App$ and there is a $\Pp\restrictedto \AA^{q^*_1}$-name 
$\name B_1\supseteq\name B_0$ of an appropriate branch for $q^*_1$.  Now
apply part (2) to $q_0, q^*_1, q_2$.

4. As in [Sh326, 3.9(2)], by induction on  the order type of
$(\AA^{q_2}\setminus\AA^{q_1})$, using (3).

5, 6.  Since (6) includes (5), it suffices to prove (6); but as we go
through the details we will treat the cases corresponding to (5) first.
We point out at the outset that if $\delta_2$ is a successor ordinal
or a limit of uncountable cofinality, then we can replace the $q_j$ by
their union, which we call $q$, 
setting $\beta=\sup_j \beta_j$, so all these cases can be
treated using the notation of (5).

We will prove by induction on $\gamma < \omega_2$ that if all
$\beta_j
\le \gamma$ and all $p_i$ belong to  $\App\restrictedto \gamma$, 
then the claim (6) holds for some $r$ in $\App\restrictedto
\gamma$.  

We first dispose of most of the special cases which fall under clause (5).
If $\delta_1 = \delta_0+1$ is a successor ordinal it suffices
to apply (4) to $p_{\delta_0}$ and $q$. 
So we assume for the present that $\delta_1$ is a limit
ordinal.  In addition if $\gamma = \beta$ we take $r=q$, so we will
assume $\beta < \gamma$ throughout.

\noindent {\sl The case $\gamma = \gamma_0+1$, a successor.}

In this case our induction hypothesis applies to the $p_i\restrictedto
\gamma_0$, $q$, $\beta$, and $\gamma_0$, yielding $r_0$ in
$\App\restrictedto \gamma_0$ with $p_i\restrictedto \gamma_0\le
r_0$ and $q\leend r_0$.  
What remains to be done is an amalgamation of $r_0$ with all of the
$p_i$, where 
$\dom p_i\includedin \dom r_0 \union \{\gamma_0\}$, 
and where one may as well
suppose that $\gamma_0$ is in $\dom p_i$ for all $i$.  This is a slight
variation on 1.9 (1 or 3) (depending on the value of
$\epsilon^{p_i}_\gamma$, which is independent of $i$). 

\noindent {\sl The case $\gamma$ a limit of cofinality greater than
$\aleph_1$.}

Since $\delta_1 < \aleph_2$  there is some $\gamma_0 < \gamma$ such 
that all $p_i$ lie in $\App\restrictedto \gamma_0$ and $\beta <
\gamma_0$, and the induction hypothesis then yields the claim.

\noindent {\sl The case $\gamma$ a limit of cofinality $\aleph_1$.}

Choose $\gamma_j$ a strictly increasing and continuous 
sequence of length at most $\omega_1$ with supremum
$\gamma$, starting with $\gamma_0= \beta$. 
By induction choose $r_j\in \App\restrictedto \gamma_j$ for
$i<\omega_1$ such that: 
$$
r_0=q;\leqno(0)
$$
$$
r_j \leend r_{j'} \hbox{ for $j < j' < \omega_1$};\leqno(1)
$$
$$
p_i\restrictedto \gamma_j\le r_j \hbox{ for $i < \delta_1$ and 
$j < \omega_1$}.\leqno (2)
$$
At successor stages the inductive hypothesis is applied to $p_i\restrictedto
\gamma_{j+1}$, $r_j$, $\gamma_j$, and $\gamma_{j+1}$.  At limit stages $j$
we apply the inductive hypothesis to $p_i\restrictedto \gamma_j$,
$r_{j'}$ for $j'<j$, $\gamma_{j'}$ for $j'<j$, and $\gamma_j$; and
here $(6)$ is used, inductively.

Finally let $r = \Union r_j$.

We now make an observation about the case of (5) that we have not yet
treated, in which $\gamma$ has cofinality $\omega$. In this case we
can use the same construction used when $\gamma$ has cofinality
$\aleph_1$, except for the last step (where we set $r=\Union r_j$,
above).  What is needed at this stage would be an instance of (6),
with the $r_j$ in the role of the $q_j$ and $\delta_2=\omega$.

This completes the induction for the cases that fall under the
notation of $(5)$, apart from the case in which $\gamma$ has
cofinality $\omega$, which we reduced to an instance of (6) with the
same value of $\gamma$ and with $\delta_2=\omega$.
Accordingly as we deal with the remaining cases
we may assume $\delta_2=\omega$.
In this case $q=\Union q_j$ is a
well-defined object, but not necessarily in \App, as the filter
$\name \FF^q$ is not necessarily an ultrafilter (there are reals
generated by $\Pp\restrictedto (\dom q)$ which do not come from any
$\Pp\restrictedto (\dom q_j)$). 

We distinguish two cases. 
If $\beta := \sup \beta_j$ is less
than $\gamma$, then induction applies, delivering an element $r_0\in
\App\restrictedto \beta$ with $p_i\restrictedto \beta \le r_0$ and all 
$q_j\leend r_0$. This $r_0$ may then play the role of $q$ in an
application of 1.9 (5).

In some sense the main case (at least as far as the failure of
continuity is concerned) is the remaining one in which
$\beta=\gamma$. Notice in this case that although $p_i\restrictedto
\beta_j\le q_j$ it does not follow that $p_i\restrictedto \beta\le q$
(for the reason mentioned above: $p_i\restrictedto \beta$ includes an
ultrafilter on part of the universe, while the filter associated with
$q$ need not be an ultrafilter).  All that is needed at this stage is
an ultrafilter containing all $\name \FF^{p_i}\union \name \FF^{q_j}$.  As
this is a directed system of filters, it suffices to check the
compatibility of each such pair, as was done in 1.9 (2).
\qed

\\ Construction, first version.

We force with \App and the generic object gives us the name of an
ultrafilter in $V[\App][\Pp]$. The forcing is $\aleph_2$-complete by 
1.9 (5).  We also claim that it satisfies the $\aleph_3$-chain condition,
and hence does not collapse cardinals and does not affect our
assumptions on cardinal arithmetic. (Subsets of $\aleph_2$ are added,
but not very many.)  In particular $(\prod ^{\AA^r}(T^1_n,T^2_n)
/\name \FF^r)^{V[\Pp\restrictedto \AA^r]}$ is a $\Pp\restrictedto
\AA^r$-name, not dependent on forcing with \App. 

We now check the chain condition. Suppose we have an antichain
$\{q_\alpha\}$ of cardinality $\aleph_3$ in \App, where for
convenience the index $\alpha$ is taken to vary over ordinals of
cofinality $\aleph_2$.  We claim that by Fodor's lemma, we may suppose that
the condition $q_\alpha\restrictedto \alpha$ is constant.  One
application of Fodor's lemma allows us to assume that $\gamma=\sup
(\AA^{q_\alpha}\intersect \alpha)$ is constant. Once $\gamma$ is
fixed, there are only $\aleph_2$ possibilities for
$q_\alpha\restrictedto \gamma$, by our assumptions on the ground
model, and a second application of Fodor's lemma allows us to take
$q_\alpha\restrictedto \gamma$ to be constant.

Now fix $\alpha_1$ of cofinality $\aleph_2$ (or more accurately, in the
set of indices which survive two applications of Fodor's lemma), 
and let $q_1=q_{\alpha_1}$, $\beta = \sup \AA^{q_1}$, and take $\alpha_2>
\beta$ of cofinality $\aleph_2$. We find that $q_2=: q_{\alpha_2}$ and
$q_1$ are compatible, by 1.9 (4), and this is a contradiction.

\\ Construction, second version.

If we wish  to apply the method of [ShHL162] (over a suitable ground model) and
build the name of our ultrafilter in the ground model, we proceed as follows.
For $\alpha\le \aleph_3$ we choose  $G^\alpha\includedin
\App\restrictedto \alpha$, directed under $\le$, inductively
as in [Sh326, \S3], making all the commitments we can;
more specifically, take $\NN\prec(H(\bet_{\omega+1}^+), \in)$ of
cardinality $\aleph_2$
with $\delta \in\NN$, $\aleph_2\subseteq \NN$,
$\NN$ is $(<\aleph_2)$-complete, and the oracle associated with
$\diamond_S$ belongs to $\NN$, and make all the commitments known to
$\NN$. Then $G^\alpha$ is in the ground model but behaves like a
generic object for $\App\restrictedto \alpha$ in $V[\Pp\restrictedto
\alpha]$, and in particular gives rise to a name $\name \FF^\alpha$.

The lengthy discussion in [Sh326 \S3] is useful for developing intuition.
Here we will just note briefly that what is called a commitment
here is really an isomorphism type of commitment, in a more
conventional sense; this is a device for compressing $\aleph_3$
possible commitments into a set of size $\aleph_2$. 

The formalism is documented in the appendix to this paper, but as we
have said it has to be adapted to allow weaker continuity axioms.
Compare paragraphs A1 and A6 of the appendix.  The axioms in the
appendix have been given in a form suitable to their application to
the proof of the relevant combinatorial theorem, rather than in the
form most convenient for verification.  1.9 above represents the sort
of formulation we use when we are actually verifying the axioms. 

We will now add a few details connecting 1.9 with the eight axioms of
paragraph A6.  The first three of these are formal and it may be
expected that they will be visibly true of any situation in which this
method would be applied. The fourth axiom is the so-called
amalgamation axiom which has been given in a slightly more detailed
form in 1.9 (4). The last four axioms are various continuity axioms,
which are instances of 1.9 (5). We reproduce them here:

\item {$5'$.}
If $(p_i)_{i < \delta} $ is an increasing sequence in \App
of length less than
$\lambda$, then it has an upper bound $q$.

\item {$6'$.}
If $(p_i)_{i < \delta}$ is an increasing sequence of length less than
$\lambda$ of members of $\App\restrictedto (\beta+1)$, with $\beta <\lambda^+$ and
if $q\in \App\restrictedto \beta$ 
satisfies $p_i\restrictedto \beta\le q$ for all
$i< \delta$, then $\{p_i:i < \delta\}\union \{q\}$ has an upper bound
$r$ in $\App$ with $q\leend r$.

\item {$7'$.}
If $(\beta_j)_{j < \delta}$ is a strictly increasing sequence of
length less than
$\lambda$, with each $\beta_j<\lambda^+$, and $p\in \App$,
$q_i\in \App\restrictedto \beta_i$, 
with $p\restrictedto \beta_j \le q_j$, and $p_{j'}\restrictedto
\beta_j = p_j$ for $j<j' < \delta$,
then $\{p\}\union \{q_j:i< \delta\}$ has an
upper bound $r$ with all $q_j\leend r$.

\item {$8'$.}
Suppose $\delta_1, \delta_2$ are limit ordinals less than $\lambda$, and
$(\beta_j)_{j< \delta_2}$ is a strictly increasing continuous sequence of
ordinals less than $\lambda^+$. Let 
$I(\delta_1,\delta_2) := (\delta_1+1)\times (\delta_2 +1)
-\{(\delta_1,\delta_2)\}$.  
Suppose that for $(i,j)\in I(\delta_1,\delta_2)$
we have $p_{ij}\in
\App\restrictedto {\beta_i}$ such that 
$$
\eqalign{
i\le i' &\implies p_{ij}\le p_{i'j}\cr
j\le j' &\implies p_{ij} = p_{ij'}\restrictedto \beta_j;\cr
}
$$
Then $\{p_{ij}: (i,j)\in I(\delta_1,\delta_2)\}$ has an upper bound
$r$ in $\App$ with $r\restrictedto \beta_j =p_{\delta_1,j}$ for all
$j<\delta_2$.

The first three are visibly instances of 1.9 (5). In the case of axiom $(8')$ 
we set
$p_i = p_{i, \delta_2}$ for $i < \delta_1$
and $q_j= p_{\delta_1, j}$ for $j < \delta_2$.
Then $p_i \restrictedto \beta_j = p_{i, j}\le q_j$, 
so 1.9 (5) applies
and yields $(8')$.

\\ Lemma

Suppose $\delta<\aleph_3$, $\cf(\delta)=\aleph_2$, and
$\Hh^\delta \subseteq \Pp\restrictedto \delta$ is generic for
$\Pp\restrictedto \delta$.
Then in $V[\Hh^\delta]$ we have:
$$\hbox{$\prod_n {}^\delta (T_n^1,T_n^2)/\name \FF^\delta[\Hh^\delta]$ is
 $\aleph_2$-compact.}$$

\proof:
 Similar to 1.8 (2).  
We can use some $\name
x_\beta$ with $\beta$ of cofinality less than $\aleph_2$ to realize
each type.  In the forcing version, this means \App forces our claim
to hold since it can't force the opposite.  In the alternative
approach, 
what we are saying is that the commitments we have made
include commitments to make our claim true. As $2^{\aleph_1}=\aleph_2$ in
$V[\Hh^\delta]$  we
can ``schedule'' the commitments conveniently, so that each particular
type of cardinality $\aleph_1$ that needs to be considered by stage
$\delta$ in fact appears before stage $\delta$. 
\qed

\\ Killing isomorphisms

We begin the verification that our filter $\name \FF$ satisfies the
condition of Proposition B. We suppose therefore that we have a
$\Pp$-name $\name F$ and a condition
$\pp^*\in \Pp$ forcing:
$$\hbox{``$\name F $ is  a map from 
$\prod_n  T^1_n$  onto $\prod_n  T_n^2$ which
represents an isomorphism modulo $\name \FF$.''}$$
We then have a stationary set $S$ of ordinals 
$\delta<\aleph_3$ of cofinality $\aleph_2$
which satisfy:
$$\pp^*\in \Pp\restrictedto \delta.
\leqno(a)$$
$$\hbox{For every $\Pp\restrictedto \delta$-name $\name x$ for an
element of $\prod_n T^1_n$,
$\name F(\name x)$ is a $\Pp\restrictedto \delta$-name.}
\leqno(b)
$$
$$\hbox{Similarly for $\name F^{-1}$.\hbox to 3.15 in{}}
\leqno(c)
$$
If we are using our second approach, over an L-like ground model:
$$\hbox{At stage $\delta$ of the construction of the 
$G^\alpha$, the diamond ``guessed'' $\pp^\delta=\pp^*$ and
$\name F_\delta=\name F \restrictedto \delta$.}
\leqno(d)$$

\noindent (In this connection, recall that 
the guesses made by diamond influence
the choice of ``commitments'' made in the construction of the
$G^\delta$.) 
Let 
$\name y^* =: \name F(\name x_\delta)$.
Then:
$$
\pp^*\forces\left\{ 
\eqalign 
{&\hbox{``$ \name y^*$ induces a branch
in $(\prod_n  T_n^2/\name \FF)^{V[\Pp\restrictedto \delta]}$ which is
the image under    $\name F_\delta$ of} \cr 
&\hbox{the branch which
$\name x_\delta$ induces on ($\prod_n T^1_n/\name \FF)^{\Vof \delta}$.''}\cr} 
\right.
\leqno (*)_{\name y^*}
$$

Now we come to one of the main points. 
We claim that there is some $q^*\in G$ with the following property:
$$
\eqalign{
&\hbox{Given $q_1\in G^\delta$ with 
$q^*\restrictedto \delta\le
q_1$ and $\Pp\restrictedto \AA^{q_1}$-names $(\name x, \name y)$ with
$\name x \in \prod T_n^1$, $\name y \in \prod T_n^2$,}\cr
&\hbox{then for any $q_3'\in \App$ with $q_1,q^*\le q_3'$ and
$q_3'\restrictedto \delta\in G^\delta$, $\pp^*$
forces:}\cr
&\hbox{``If $\name y = \name F(\name x)$ then $\name x \le \name
x_\delta$ iff $\name y \le \name y^*$, and}\cr
&\hbox{ if $\name y$ and $\name
F(\name x)$ are incomparable, then $\name x\le \name x_\delta$ implies
$\name y \not \le \name y^*$.''}\cr
}
\leqno (\dag)_\delta
$$
Notice here that $q_3'$ need not be in $G$. 

The reason for this depends  slightly on which of the two
approaches to the construction of $G$ we have taken.  In a straight
forcing approach, we may say that some $q^*\in G$ forces 
$(*)_{\name y^*}$, and this yields $(\dag)_\delta$.
In the second, pseudo-forcing, approach we find that our 
``commitments'' include a
commitment to falsify $(*)_{\name y^*}$ if possible; as 
we did not do so, at a 
certain point it must have been impossible to falsify it,
which again translates into $(\dag)_\delta$.

We now fix $q^*$
satisfying $(\dag)_\delta$, and we set $q_0=q^*\restrictedto \delta$.
At this stage, $(\dag)_\delta$ gives some sort of local definition of $\name
F\restrictedto \delta$, on a cone in $(\prod {}^\delta T_n^1/\name
\FF^\delta)^{\Vof \delta}$ determined by $q_0$. The next result allows
us to put this definition in a more useful form (and this is nailed
down in 1.15). One may think of this as an elimination of quantifiers.

\\ Lemma

\proclaim{} Suppose that: 
$$\hbox{$q_0,q_1,q_2,q_3$ are in $\App$ with
$q_0=q_2\restrictedto \beta_0\le q_1\leend q_3$, 
and $q_2\le q_3$.}\leqno(1)$$ 
$$\hbox{$q_0\le r_0\in \App$ with $\AA^{q_1}\includedin
\AA^{r_0}\includedin \beta_0$.}
\leqno(2)$$
Let $\AA_i=\AA^{q_i}$ for $i=0,1,2,3$,
and suppose that 
$$\left\{
\eqalign{&\hbox{$\name f_0$ is a 
$\Pp\restrictedto \AA^{r_0}$-name of a map from 
$(\prod_n (T^1_n,T_n^2))^{\Vof {\AA_{1}}}$
to $(\prod_n (T^1_n,T_n^2))^{\Vof {\AA^{r_0}}}$}\cr
&\hbox{ representing a partial elementary embedding of}\cr 
&\hbox{$(\prod_n{}^{\AA_0} (T^1_n,T_n^2)/\name \FF\restrictedto
\AA_1)^{\Vof {\AA_{1}}}$ into 
$(\prod_n {}^{\AA_0}(T^1_n,T_n^2)/\name
\FF\restrictedto\AA^{r_0})^{\Vof {\AA^{r_0}}}$ }\cr
&\hbox{which is equal to the identity on
$(\prod_n(T^1_n,T_n^2)/\name \FF\restrictedto \AA_0)^{\Vof
{\AA_0}}$.}\cr}\right.\leqno(3)
$$
Then there is an $r\in \App$ with:
$$q_2\le r;\quad r_0\leend r;\quad \AA_3\includedin \AA^{r};\quad
\AA^{r}\intersect \beta_0=\AA^{r_0};$$
and there is a $\Pp$-name $\name f$ of a function from 
$(\prod_n(T^1_n,T_n^2))^{\Vof {\AA_3}}$ into 
$(\prod_n(T^1_n,T_n^2))^{\Vof {\AA^r}}$
representing an elementary embedding of 
$^{\AA_2}(\prod_n(T^1_n,T_n^2)/\name \FF\restrictedto \AA_3)^{\Vof {\AA_3}}$
into 
$^{\AA_2}(\prod_n(T^1_n,T_n^2)/\name \FF\restrictedto \AA^r)^{\Vof  {\AA^r}}$
which is the identity on 
$(\prod_n(T^1_n,T_n^2)/\name \FF\restrictedto \AA_2)^{\Vof {\AA_2}}$.
\endproclaim

\proof:
It will be enough to get $\name f$ as a partial elementary embedding,
as one may then iterate 1.8 (3) $\aleph_1$ times.

We may suppose $\beta_0 = \inf\, (\AA_3-\AA^{r_0})$. 
Let $\AA_3\setminus \beta_0=(\beta_i)_{i<\xi}$ enumerated in
increasing order. We will construct two increasing sequences, one of
names $\name f_i$ and and one of elements
$r_i\in \App$, indexed by $i\le \xi$, such that our claim holds for $\name f_i,
q_2\restrictedto \beta_i, q_3\restrictedto \beta_i, r_i$, and in
addition $\AA^{r_i}\includedin \beta_i$.
At the end we take $r=r_\xi$ and $\name f=\name f_\xi$.

\smallskip
\noindent {\sl The case $i = 0$}

Initially $r_0$ and $\name f_0$ are given.

\smallskip
\noindent{\sl The limit case}

Suppose first that $i$ is a limit ordinal of cofinality $\aleph_0$,
and
let $\AA=\Union_{j<i} \AA^{r_j}$. In this case $\Union_{j<i}\name
\FF^{r_j}$ is not an ultrafilter in $V[\Pp\restrictedto \AA]$ and
the main point will be to
prove that there is a $\Pp\restrictedto \AA$-name for an ultrafilter
$\name \FF_{i}$ extending   
$\name \FF^{q_2\restrictedto \beta_i}$ and $\Union_{j<i} \name \FF^{r_i}$, such
that 
$$
\left\{
\eqalign{
&\hbox{The map $\name f_{i}$ defined as the identity on
$(\prod_n (T^1_n,T_n^2)) ^{\Vof {(\AA_2\intersect \beta_i)}}$ and as 
$\Union_{j<i} \name f_j$ }
\cr 
&\hbox{on the latter's domain is a partial elementary map from}
\cr
&\hbox{$(\prod_n ^{\AA_2\intersect \beta_i}(T_n^1,T_n^2)/\name
\FF\restrictedto (\AA_3\intersect \beta_i)) ^{\Vof
{(\AA_3\intersect\beta_i)}}$ into
$(\prod_n ^{\AA_2\intersect \beta_i}(T_n^1,T_n^2)/\name \FF^{r_i})^
{\Vof{\AA}}$}.
\cr
}\right.\leqno(*)
$$

So it will suffice to find $\name \FF_{i}$ making $(*)$ true. This
means we must check the finite intersection property for a certain
family of (names of) sets. 
Suppose toward a contradiction that we have a condition $\pp\in
\Pp\restrictedto \AA$ forcing ``$\name a\intersect \name b\intersect
\name c=\emptyset$,'' where:
$$\def\offsetsize{.36 true in}
\eqalign{
\offset{(A)}&\hbox{$\name a$ is}
\hbox{ a $\Pp\restrictedto \AA^{r_j}$-name for a member of $\name
\FF^{r_j}$} 
\cr 
\offset{(B)}\hfill&\hbox{$\name b$ is}
\hbox{ a $\Pp\restrictedto \AA^{q_2\restrictedto \beta_i}$-name for a member of 
$\name \FF^{q_2\restrictedto \beta_i}$}\cr
\offset{(C)}\hfill &\hbox{$\name c$ is}\hbox{ the name of a set of the form:
$\{n:(T^1_n,T_n^2)\satisfies 
\name \phi(\name \xbar(n),\name f_j(\name \ybar)(n))\}$.}\cr
\offset{(C1)}&\hbox{$\name \xbar$, $\name \ybar$ are finite sequences from 
$(\prod_n (T_n^1,T_n^2))^{\Vof {\AA^{q_2\restrictedto \beta_i}}}$
and $(\prod_n (T^1_n,T_n^2))^{\Vof {(\AA_3\intersect \beta_j)}}$ 
respectively.}
\cr
\offset{(C2)}&\hbox{$\name \phi$ is a $\Pp\restrictedto \AA^{q_2\restrictedto
\beta_i}$-name for a formula in the language of $\prod_n
{}^{\AA^{q_2\restrictedto \beta_i}}(T_n^1,T_n^2)$}
\cr
\offset{(C3)}&\hbox{$\name \phi(\name \xbar, \name \ybar)$
holds in 
${}^{\AA_2\intersect \beta_i}(\prod_n (T^1_n,T_n^2)/\name \FF\restrictedto
(\AA_3\intersect \beta_i))^{\Vof {\AA_3\intersect \beta_i}}$.}
\cr
} 
$$

Here $j<i$ arises as the supremum of finitely many values below $i$.
As $\name \xbar$ can be absorbed into the language, we will drop it.

Now let $\Hh$ be generic for $\Pp\restrictedto (\AA_2\intersect \beta_j)$ 
with $\pp\restrictedto (\AA_2\intersect \beta_j)\in \Hh$, and define:
$$
\name A_n=:\{\ubar:\hbox{for some $\pp_2\ge \pp\restrictedto
(\AA_2\intersect \beta_i)$ with $\pp_2\restrictedto (\AA_2\intersect \beta_j)
\in \Hh$,\quad $\pp_2\forces$ ``$n \in \name b$ and $(T^1_n,T_n^2)\satisfies
\name \phi(\ubar)$.''}\} 
$$

$\name A_n$ is a $\Pp\restrictedto (\AA_2\intersect \beta_i)$-name of
a subset of $T_n^2$. 
Take $(\name A_n)$ as a
relation in $\prod {}^{\AA^{q_2\restrictedto \beta_j}}(T_n^1,T_n^2)$. 
By hypothesis $\{n: (T_n^1,T_n^2)\satisfies \name \phi(\name
\ybar(n))\}\in \name \FF^{q_3\restrictedto \beta_i}$, and this set is
contained in the set $\name a'=\{n:\name \ybar(n)\in \name A_n\}$, which 
belongs to
$V[\Pp\restrictedto (\AA_3\intersect \beta_j)]$. Therefore $\name a'\in
\FF^{q_3\restrictedto \beta_j}$ and applying $\name f_j$, we find:
$$ 
\{n:\name f_j(\name \ybar)(n)\in \name A_n\}\in \name \FF^{r_j}.
$$
Hence we may suppose that $\pp$ forces: for $n\in \name a$, $\name
f_j(\name \ybar)(n)\in \name A_n$. But then any element of $\name a$
can be forced by an extension of $\pp$ to lie in $\name b \intersect
\name c$, by amalgamating appropriate conditions over $\AA_2\intersect \beta_j$.

Limits of larger cofinality are easier. 

\smallskip
\noindent {\sl The successor case}

Suppose now that $i=j+1$. We may suppose that $\beta_j\in \AA_2$ as
otherwise there is nothing to prove.
If $\epsilon^{q_2}_{\beta_j}=0$ we argue as in the
previous case.
So suppose that 
$\epsilon^{q_2}_{\beta_j}=1$. In particular $\beta_j$ has cofinality
$\aleph_2$. 

Using 1.8 (3) repeatedly, and the limit case, 
we can find $\name B$, $q_1',r',\name f'$ such that:
$$
q_3\restrictedto \beta_j\leend q_1';\quad \AA^{q_1'}\includedin \beta_j;
\leqno(1)
$$ 
$$
r_j\leend r';\quad \AA^{r'}\includedin \beta_j;
\leqno (2)
$$

$$
\left\{\eqalign{
&\hbox{$ \name f'$ is a map from $\prod_n (T^1_n,T_n^2)^{\Vof {\AA^{q_1'}}}$
onto $\prod_n (T^1_n,T_n^2)^{\Vof {\AA^{r'}}}$ representing an isomorphism}\cr
&\hbox{of
$(\prod_n (T^1_n,T_n^2)/\name \FF^{q_1'})^{\Vof {\AA^{q_1'}}}$ with
$(\prod_n (T^1_n,T_n^2)/\name \FF^{r'})^{\Vof {\AA^{r'}}}$
extending $f_j$;}\cr}\right.
\leqno (3)$$
$$
\hbox{$\name B$ is a $\Pp\restrictedto \AA^{q_1'}$-name of a branch
of $(\prod_n T^1_n/\name \FF^{q_1'})^{\Vof {\AA^{q_1'}}}$
which is sufficiently generic;}\leqno(4)
$$
$$
\hbox{$\name f'[\name B]$ is a $\Pp\restrictedto \AA^{r'}$-name of a branch of
$(\prod_n T^1_n/\name \FF^{r'})^{\Vof {\AA^{r'}}}$
which is sufficiently generic.}
\leqno(5)$$

Let $q_3'$ satisfy $q_3\restrictedto \beta_i\le q_3'$, $q_1'\leend
q_3'$, with $\AA^{q_3'}\includedin \beta_i$ such that according to
$q_3'$ the vertex $\name x_{\beta_j}$ lies above $\name B$ (using 1.9(2)). 
We intend to have $r_i$ 
put $\name x_{\beta_j}$ above $\name f'[\name B]$ (to meet conditions
(5.2, 5.3) in the definition of $\App$), while meeting our other responsibilities.
As usual the problem is
to verify the finite intersection property for a certain family of
names of sets. Suppose therefore toward a contradiction that we have a condition
$\pp\in \Pp$ forcing ``$\name a\intersect\name b\intersect\name
c\intersect \name d=\emptyset$,'' where 
$$\eqalign{
\hbox{$\name a$ is}&
\hbox{ a $\Pp\restrictedto \AA^{r'}$-name of a member of $\name
\FF^{r'}$;}\cr
\hbox{$\name b$ is}&
\hbox{ a $\Pp\restrictedto \AA^{q_2\restrictedto \beta_i}$-name of a
member of $\name \FF^{q_2\restrictedto \beta_i}$;}\cr
\hbox{ $\name c$ is}&
\hbox{ the name of a set of the form $\{n:(T^1_n,T_n^2)\satisfies \name \phi(\name
x_{\beta_j}(n),\name f'(\name \ybar)(n))\}$}\cr
\hbox{$\name d$ is}&
\ \{n:T^1_n \satisfies \name x(n)<\name x_{\beta_j}(n)\}\cr
}
$$
where in connection with $\name c$ we have: 
$$
\eqalign{
&\name \ybar \in (\prod_n (T^1_n,T_n^2))^{\Vof {\AA^{q_1'}}},\cr
&\hbox{$\name \phi(\name x_{\beta_j},\name \ybar)$ is defined and
holds in  
$(\prod_n {}^{\AA_2\intersect \beta_i} (T^1_n,T_n^2)/\name
\FF^{q_3'})^{\Vof{\AA^{q_3'}}}$,}\cr} 
$$ 
and we have absorbed some parameters occurring in $\name \phi$ into
the expanded language which is associated with $V[\Pp\restrictedto
(\AA_2\intersect \beta_j)]$ as individual constants, 
while in connection with $\name d$ we have:
$$
\hbox{$\name x$ is a $\Pp\restrictedto \AA^{q_1'}$-name for a member
of $\name f'[\name B]$. 
}
$$

Let $\Hh^*\includedin \Pp$ be generic over $V$ with $\Hh\includedin
\Hh^*$ and $\pp\in \Hh^*$. Set 
$\Hh=\Hh^* \restrictedto \AA^{q_2\restrictedto \beta_j}$, 
$\Hh_1=\Hh^* \restrictedto \AA^{q_1'}$, and
$\Hh_3=\Hh^* \restrictedto \AA^{q_3'}$.
In $V[\Hh]$ we define:
$$
\eqalign
{\name A_n^1=:\,&
\hbox{$\{(x, \ubar):$ 
For some $\pp_1\in \Pp\restrictedto \AA^{r'}$, with $\pp_1\ge
\pp\restrictedto \AA^{r'}$ and   
$\pp_1\restrictedto \AA^{q_2\restrictedto \beta_j}\in \Hh$,}\cr
&\qquad\qquad\hbox{$\pp_1$ forces: ``$n\in \name a$, $\name x(n)=x$,
$\name f'(\name \ybar)(n) = \ubar$.$\}$}\cr
\name A_n^2=:\,&
\hbox{$\{(x^*,\ubar):$ 
For some $\pp_2\in \Pp\restrictedto (\AA_2\intersect \beta_i)$
with $\pp_2\ge\pp\restrictedto (\AA_2\intersect \beta_i)$ and
$\pp_2\restrictedto (\AA_2\intersect \beta_j)\in\Hh$,}\cr 
&\qquad\qquad\hbox {$\pp_2$ forces: 
``$n\in \name b$, $\name x_{\beta_j}(n)=x^*$, and 
$\name \phi(x^*,\ubar)$.''$\}$}\cr}  
$$

In $V[\Hh]$ there is no $n$ satisfying:
$$\hbox{ $\exists x, x^*, \ubar \quad (x, \ubar) \in \name
A^1_n \and (x^*,\ubar) \in \name A_n^2 \and x< x^*$.}
\leqno(*)$$
Otherwise we could extend $\pp$ by amalgamating suitable conditions
$\pp_1$, $\pp_2$, to force such an $n$ into $\name a\intersect \name b
\intersect \name c \intersect \name d$.

For $n<\omega$ and $u\in T_n^1$ let
$$
\eqalign{
\name A^2_n(\ubar) =:&
\{x\in T_n^1:(x,\ubar)\in A_n^2\}\cr
\name A^3_n(\ubar)=:&
\{x\in T_n^1: \hbox{Either $(x,\ubar)\in \name A_n^2$ or 
there is no $x'$ above $x$ in $T^1_n$ for which $(x',\ubar)\in \name
A_n^2$}\}\cr
}
$$
Then $\name A^3_n(\ubar)$ is dense in $T^1_n$.
and hence so is $\name A^3=:\prod \name A^3_n/\name
\FF^{q_2\restrictedto \beta_i}[\Hh]$. 

Let $\TT = (T^1,T^2;A^2,A^3)$ be the ultraproduct $(\prod_n
(T_n^1,T_n^2;\name A_n^2,\name A_n^3)/\FF^{q_1'})^{V[\Hh_1]}$.  
Now $\name \phi[\name x_\beta, \name \ybar]$ holds in
$\prod{}^{\AA_2\intersect \beta_i}(T_n^1,T_n^2)/\FF^{q_3'}[\Hh_3]$, so
$\name x_\beta[\Hh_3]\in A^2(\name \ybar[\Hh_3])$ (using \L o\'s'
theorem to keep track of the meaning of $A^2$ in this model). 
By the choice of $\name B$, $\name B[\Hh_1]$ meets $A^3(\name \ybar[\Hh_1])$ and
indeed:
$$
\hbox{$A^3(\name \ybar[\Hh_1]) \intersect \name B[\Hh_1]$ is unbounded in $\name
B[\Hh_1]$}
\leqno(1)
$$
For $\name z\in A^3(\name\ybar[\Hh_1])\intersect \name B[\Hh_1]$, 
as $\name z < \name x_{\beta_j}$ we have 
also $\name z\in A^2(\name \ybar[\Hh_1])\intersect \name B[\Hh_1]$. Hence in
$V[\Hh_1]$ we have:
$$
\hbox{$A^2[\name \ybar]\intersect \name B[\Hh_1]$ is unbounded in
$\name B[\Hh_1]$}\leqno(2)
$$
and hence $A^2(\name f'(\name \ybar))\intersect \name f'[\name
B][\Hh^*\restrictedto \AA^{r'}]$ is unbounded in $\name f'[\name
B][\Hh^*\restrictedto \AA^{r'}]$, and  we can find $\name z\in
A^2(\name f'(\name \ybar[\Hh_3]))\intersect \name f'[\name
B][\Hh^*\restrictedto \AA^{r'}]$ with $\name x < \name z$ in $\prod_n
T_n^1/\FF^{r'}[\Hh^*\restrictedto \AA^{r'}]$.

In particular for some $n\in \name a[\Hh^*]$, we have $\name x(n)[\Hh^*] <
\name z(n)[\Hh^*]$ in 
$T_n^1$ and $\name z(n)\in A^2(\name \ybar(n))$. Letting
$x=\name x(n)[\Hh_1]$, $x^*=\name z(n)[\Hh_1]$, and $u=\name f'(\name
\ybar)(n)[\Hh\restrictedto \AA^{r'}]$,  we find that $(*)$ 
holds in $V[\Hh]$, a contradiction. 
\qed

\\ Weak definability

\proclaim{Proposition}
Let $\delta < \aleph_3$ be an ordinal of cofinality $\aleph_2$
satisfying conditions 1.13 (a-d). 
Suppose $q_1,q_2\in G$,  $q_2 \restrictedto \delta=q_0\le q_1$,
$\AA^{q_1}\subseteq\delta$, $\delta \in A^{q_2}$, 
$\name y^*$ is a $\Pp\restrictedto \AA^{q_2}$-name of an element of
$\prod_{n}T_n^2$,  
and $\epsilon_\delta^{q_2}=1$. 
Suppose further that $\name x'$, $\name x''$ and $\name y'$,
$\name y'' $
are $\Pp\restrictedto \AA^{q_1}$-names,
 $\pp\in \Pp$, $\pp_i=\pp\restrictedto \AA^{q_i}$ $(i=1,2)$, and:
$$
\eqalign
{\pp_1\forces&
\hbox{``$\name x',\name x''\in \prod_n T^1_n$, and 
$\name y', \name y''\in \prod_n T_n^2$;''\quad $\pp_2\forces$ ``$\name F(\name
x_\delta) = \name y^*$''}\cr 
\pp_1\forces&\hbox{``The types of $(\name x', \name y')$ and of $(\name
x'', \name y'')$ over $\{\name x /\name \FF:\name x $ a
$\Pp\restrictedto \AA^{q_0}$-name of a 
member of\/ $\prod_n {}^{\AA^{q_0}}(T^1_n, T_n^2)\}$}\cr
&\hbox{in the model $(\prod_n {}^{\AA^{q_0}}(T^0_n, T^1_n)/\name
\FF^{q_1})^{\Vof {\AA^{q_1}}}$ are equal.''}\cr}
$$

Then the following are equivalent.
\item{1.} There is $r^0\in \App$ such that $q_1,q_2\le r^0$,
 $r^0\restrictedto \delta\in G^\delta$, and  
$$\pp\forces
\hbox{``$\prod_n T^1_n/\name \FF^{r^0}\models 
(\name x'/\name \FF^{r^0} <\name x_\delta/\name \FF^{r^0})$ 
and $\prod_n T_n^2/\name \FF^{r^0}\models
(\name y'/\name \FF^{r^0}<\name y^*/\name \FF^{r^0})$'';}$$ 
\item{2.} There is $r^1\in \App$ such that $q_1,q_2\le r^1$,
$r^1\restrictedto \delta\in G^\delta$
and  
$$\pp\forces
\hbox{`` $\prod_n T^1_n/\name \FF^{r^1}\models
(\name x''/\name \FF^{r^1}<\name x_\delta/\name \FF^{r^1})$
and $\prod_n T_n^2/\name \FF^{r^1}\models
(\name y''/\name \FF^{r^1}<\name y^*/\name \FF^{r^1})$.''}$$
\endproclaim

\proof:
It suffices to show that (1) implies (2). 
Take $\Hh^\delta \subseteq \Pp\restrictedto \delta$ generic over $V$
with $\pp_1\in \Hh^\delta$, and suppose that $r^0$ is as in (1). 
Let $r_0=r^0\restrictedto \delta$ and let $\name f_0$ be the extension
of the identity
map on $(\prod T_n^1)^{\Vof \AA^{q_0}}$ by: 
$\name f_0(\name x') = \name x''$, 
$\name f_0(\name y') = \name y''$.
Writing $\beta_0=\delta$ and taking $q_3$ provided by 1.9 (4), 
we recover the assumptions of 1.13, which
produces a certain $r$ in $\App$, an end extension of $r_0$;
here we may easily keep $r\restrictedto \delta\in G^\delta$ (cf. 1.12). 
It suffices to take $r^1=r$.
\qed

\\ Definability.

We claim now that $\name F$ is definable on a cone
by a first order formula.
For a stationary set $S_0$ of $\delta<\aleph_3$ of cofinality $\aleph_2$, we
will have conditions (a-d) of 1.13 which may be expressed as follows:
$$\hbox{Both
$\name F\restrictedto (\Pp\restrictedto \delta-\hbox{names})$
and $\name F^{-1}\restrictedto (\Pp\restrictedto \delta-\hbox{names})$
are $\Pp\restrictedto \delta$-names;}$$
When working with $\diamond_S$:
$$
\hbox{$\diamond_S$ 
guessed the names of these two restrictions and also guessed $\pp^*$ 
correctly;}$$
and hence for suitable $\name y_\delta$ and $q_\delta^*$ we have the
corresponding conditions $(*)_{\name y_\delta}$ and $(\dag)_\delta$
(with $q_\delta^*$ in place of $q^*$). By Fodor's lemma, on a stationary set
$S_1\includedin S_0$ we have $q_0=q_\delta^*\restrictedto \delta$
is constant, and also the isomorphism type of the pair
$(q_\delta^*,\name y_\delta)$ over $\AA^{q_0}$ is constant.

So 
for $\delta$ in $S_1$, we have the following two properties, holding
for $\name x'$ in $V[\Pp\restrictedto \delta]$ and $\name y' = \name
F(\name x')$),   
by $(\dag)_\delta$ and 
1.15 respectively:
$$\eqalign{
&\hbox{1. The decision to put $\name x'$ 
below $\name x_\delta$ implies
also that $\name y'$ must be put below $\name y^*$; and}\cr
&\hbox{2. This decision is determined by the type of $(\name x',\name
y')$ in $\prod {}^{\AA^{q_0}}(T^1_n,T_n^2)/\name \FF^{V[\Hh][\Pp\restrictedto
\delta/\Hh]}$.}\cr}$$
As $S_1$ is unbounded below $\aleph_3$ this holds generally.

This gives a definition by types of the isomorphism 
$\name F$ above
the branch in $\prod T^1_n/\name \FF^{\Vof{\AA^{q_0}}}$ which the
condition $q_\delta^*$ 
says that the vertex $\name x_\delta$ induces there (using 1.9 (2)),
and this branch does not 
depend on $\delta$. Note that this set contains a cone, and the
image of this cone is a cone in the image.
Now by
$\aleph_2$-saturation of
$\prod_n {}^{\AA^{q_0}}(T^1_n,T_n^2)/\name \FF^{\Vof \AA}$
we get a first order definition on a smaller cone; 
this last step is written out in
detail in the next paragraph. 
This proves Proposition B.

\\ Lemma {\sl (true definability)}

\proclaim{} Let $M$ be a $\lambda$-saturated structure,
and $A\includedin M$ with $|A|<\lambda$.
Let $(D_1;<_1)$, $(D_2; <_2)$ be $A$-definable trees in $M$; that is,
the partial orderings $<_i$ are linear below each node. Assume that
every node of $D_1$ or $D_2$ has at least two immediate successors. 
Let $F:D_1\to D_2$ be a tree isomorphism which is type-definable in
the following sense:
$$
\left [f(x)=y \and \tp(x,y/A)=\tp(x',y'/A)\right ]\implies f(x')=y'.
$$
Then $f$ is $A$-definable, on some cone of $D_1$.
\endproclaim

Before entering into the proof, we note that we use somewhat less
information about $F$ (and its domain and range) than is actually
assumed; and this would be useful in working out the most general form
of results of this type (which will apply to some extent in any
unsuperstable situation). We intend to develop this further elsewhere,
as it would be too cumbersome for our present purpose.

The proof may be summarized as follows. If a function $F$ is definable by
types in a somewhat saturated model, then on the locus of each 1-type,
it agrees with the restriction of a definable function. If $F$
is an automorphism and the locus of some 1-type separates the points in a
definable set $C$ in an appropriate sense, then $F$ can be recovered,
definably, on $C$. Finally, in sufficiently saturated trees of the
type under consideration, some 1-type separates the points of a cone.
Details follow.

\proof:
If we replace $M$ by a $\lambda$-saturated
elementary extension, the definition of $F$
by types continues to work (and the extension is an elementary
extension for the expansion by $F$). In particular, replacing $|M|$ by
a more saturated structure, if necessary, but keeping $A$ fixed, 
we may suppose that $\lambda>|T|,|A|,\aleph_0$. 

We show first:
$$\hbox{There is a 1-type $p$ defined over $A$ such that
its set of realizations 
$p[D_1]$ is dense in a cone of $D_1$,}\leqno (1)
$$
i.e., for some $a$ in $D_1$ we require that any element above $a$ lies
below a realization of $p$. For any 1-type $p$ over $A$, if
$p[D_1]$ does not contain a cone of $D_1$ then by saturation there is
some $\phi\in p$ with:
$$
\forall a\exists b>a\, \neg\exists x >b\, \phi(x)
$$
So if $(1)$ fails we may choose one such formula $\phi_p$ for each
1-type $p$ over $A$, and then it is consistent (hence true) that we have a
wellordered increasing sequence $a_p$ (in the tree ordering) such that
for each $1$-type $p$, above $a_p$ we have:
$$
\neg\exists x>a_p\, \phi_p(x)
$$
By saturation there is a further element $a$ above all $a_p$ (either by
increasing $\lambda$ or by paying attention to what we are actually doing)
and we have arranged that there is no 1-type left for it to realize.
As this is improbable, (1) holds. We fix a 1-type $p$ and an element
$a_0$ in $D_1$ so that the realizations of $p$ are dense in the cone
above $a_0$.  It is important to note at this point that the density
implies that any two distinct vertices above $a_0$ are separated by
the realizations of $p$ in the sense that there is a realization of
$p$ lying
above one but not the other (here we use the immediate splitting
condition we have assumed in the tree $D_1$).

Let $a$ realize the type $p$, and let $q$ be the type of $a,F(a)$
over $A$. If $b$ is any other realization of $p$, then there is an 
element $c$ with $b,c$ realizing $q$, and hence $F(b)=c$; thus $p$
determines $q$ uniquely. Furthermore each realization $a$ of $p$
determines a unique element $b$ such that $a,b$ realizes $q$, and
hence by saturation there is a formula $\phi(x,y)\in q$ so that
$\phi(x,y)\implies \exists!z \, \phi(x,z)$. Hence $p\union \{\phi\}\proves
q$.

Now the following holds in $M$:
$$p(x)\union p(x') \union \{\phi(x,y),\phi(x',y')\}\implies (x<x'\iff
y<y')$$
and hence for some formula $\alpha(x)\in p$ the same holds with $p$
replaced by $\alpha$. We may suppose $\phi(x,y)\implies \alpha(x)$ and
conclude that $\phi(x,y)$ defines a partial isomorphism $f$. Let $B$ be
$\{a>a_0:\exists y \phi(a,y)\}$.
$f$ coincides with $F$ on the set of realizations of $p$ above
$a$, and the action of $F$ on this set determines its action on the
cone above $a$ by density (or really by the separation condition
mentioned above), so $f$ coincides with $F$ on $B$. Furthermore the
action of $F$ on $B$ determines its action on the cone above $a_0$
definably, so $F$ is definable above $a$.

The definition $\phi^*(x,y)$ of $F$ on the cone above $a$ 
obtained in this manner may easily be written down
explicitly:
$$
\hbox{``}
\forall x', y'\, 
\left [\phi(x',y')\implies (x<x'\iff y < y')\right ]
\hbox{''}
$$
\qed

For the application in 1.16 we take $\lambda = \aleph_2$. 

\\ Remark

\proclaim {Proposition} 
$\Pp$ forces: In $\prod_n T^1_n/\name \FF$ ($\name \FF=\name
\FF[G^{\aleph_3}]$), every full branch is an ultraproduct of branches
in the original trees $T^1_n$. 
\endproclaim
\proof (in brief):
One can follow the line of the previous argument, or derive the result
from Proposition B. Following the line of the previous argument we
argue as follows.
If $\name B$ is a $\Pp$-name for such a branch, then for a stationary
set of ordinals 
$\delta<\aleph_3$ of cofinality $\aleph_2$, $\name B \intersect
(\prod_n T^1_n/\name \FF)^{\Vof\delta}$ will be a full
branch and a $\Pp\restrictedto \delta$-name, guessed correctly by
$\diamond_S$. We tried to make a commitment to terminate this
branch, but failed, and hence for some $q^*$ and $y^*$ witnesses to
the failure, we were unable to omit having $q^*\restrictedto \delta\in
G^\delta$ where $q^*$ is essentially the support of ``$y^*$ is a
bound''. Using 1.14 one shows that the branch was definable at this
point by types in $\aleph_1$ parameters, and by $\aleph_2$-compactness
we get a first order definition, which by Fodor's lemma can be made
independent of $\delta$.
\qed

Filling in the details in the foregoing argument constitutes an
excellent, morally uplifting exercise for the reader. However the more
pragmatic reader may prefer the following derivation of the
proposition from Proposition B. 

In the first place, we may replace the trees $T_n^1$ in the
proposition above by the universal tree of this type, which we take to
be $T=\Zz^{<\omega}$ (writing $\Zz$ rather than $\omega$ for the sake
of the notation used below).
Now apply Proposition B to the pair of sequences $(T_n^1)$, $(T_n^2)$
in which $T_n^i=T$ for all $i,n$. 
Using the model of ZFC and the ultrafilter 
referred to in Proposition B, suppose $B$ is a full branch of $T^*=\prod
T_n^2/\FF$, and let $\Zz^*=\Zz^\omega/\FF$, $\Nn^*=\Nn^\omega/\FF$.
For each $i\in\Nn^*$ let $B_i$ be the $i$-th node of $B$;
this is a sequence in $(\Zz^*)^{[0,i]}$ which is coded in $\Nn^*$. 
Define an automorphism $f_B$ of $T^*$ whose action on the $i$-th level
is via addition of $B_i$ (pointwise addition of sequences). Applying
Proposition B and \L o\'s' theorem to this automorphism, 
we see that $f_B$ is the ultraproduct of addition maps corresponding
to various branches of $T$, and that $B$ is the ultraproduct of these
branches. 

\\ Corollary

\proclaim{} It is consistent with ZFC that
$\Rr^\omega/\FF$ is Scott-complete for some ultrafilter $\FF$.
\endproclaim

Here $\Rr^\omega/\FF$ is called {\sl Scott-complete} if it has no
proper dedekind cut $(A,B)$ in which $\inf(b-a:a\in A,b\in B)$ is $0$
in $\Rr^\omega/\FF$. 
1.18 is sufficient for this by [KeSc, Prop. 1.3]. This corollary
answers Question 4.3 of [KeSc, p. 1024].

\\ \proclaim{Remark} \endproclaim

The predicate ``at the same level'' may be omitted from the language
of the trees $T_n^i$ throughout as the condition on $\name x_\delta$
that uses this (the ``full branch'' condition) follows from the
``bigness'' condition: meeting every suitable dense subset.


GARBAGE HEAP:
From 1.9.

\item{5.} Assume $\delta < \aleph_2$, that $(q_i)_{i< \delta}$ is an
increasing sequence from \App, that $(\beta_i)_{i<\delta}$ is a strictly
increasing sequence of ordinals, and that $(p_i)_{i<\delta}$ satisfies:
$$
\hbox 
   {For $i<\delta$: $q_i\restrictedto \beta\le p_i\in
    \App\restrictedto \beta_i$;%
   }
\qquad
\hbox
   {For $i<j<\delta$: $p_i\leend p_j$.%
   }
$$
Then there is an 
$r\in \App$    with 
$p_i \leend r$ and $q_i\le r$    for all    $i< \delta$. 
If each   $q_i$   belongs to  $\App_{\sup \beta_i}$  then  $r$  may be taken to
have domain  $\Union_i(\dom q_i\union \dom p_i)$.

5. We will prove by induction on $\gamma < \omega_2$ that if $p_i,
q_i\in \App\restrictedto \gamma$ and for all $i$ we have $\beta_i\le
\gamma$, then the claim holds (with $r$ in $\App\restrictedto
\gamma$).  If $\delta = \delta_0+1$ is a successor ordinal it suffices
to apply (4) to $q_{\delta_0}$ and $p_{\delta_0}$, with $\beta =
\beta_{\delta_0}$. So we assume throughout that $\delta$ is a limit
ordinal. In particular $\beta_i < \gamma$ for all $i$.

\noindent {\sl The case $\gamma = \gamma_0+1$, a successor.}

In this case our induction hypothesis applies to the $q_i\restrictedto
\gamma_0$, the $p_i$, the $\beta_i$, and $\gamma_0$, yielding $r_0$ in
$\App\restrictedto \gamma_0$ with $p_i, q_i\restrictedto \gamma_0\le
r_0$ (and with a side condition on the domain if all $q_i\restrictedto
\gamma_0$ lie in $\App\restrictedto (\sup \beta_i)$).  What remains
then is an amalgamation of $r_0$ with all of the $q_i$, where $\dom
q_i\includedin \dom r \union \{\gamma_0\}$, and where one may as well
suppose that $\gamma_0$ is in $\dom q_i$ for all $i$.  This is a slight
variation on 1.9 (2,3) (depending on the value of
$\epsilon^{q_i}_\gamma$, which is independent of $i$). 

\noindent {\sl The case $\gamma$ a limit of cofinality greater than
$\aleph_1$.}

Since $\delta < \aleph_2$  there is some $\gamma_0 < \gamma$ such 
that all $p_i, q_i\in \App\restrictedto \gamma_0$ and all $\beta_i <
\gamma_0$, and the induction hypothesis then yields the claim.

\noindent {\sl The case $\gamma$ a limit of cofinality $\aleph_1$.}

If $\gamma = \sup \beta_i$ then $r=\Union p_i$ suffices. Assume
therefore that $\gamma_0 := \sup \beta_i < \gamma$.  By the induction
hypothesis applied to $q_i\restrictedto \beta_i$, $p_i$, and
$\gamma_0$, we have $r_0\in \App\restrictedto \gamma_0$ with
$q_i\restrictedto \gamma_0, p_i \le r_0$ and $\dom r_0 = \Union_i
(\dom q_i\restrictedto \gamma_0 \union \dom p_i)$. 

Choose $\gamma_i^*$ a strictly increasing and continuous 
sequence of length $\omega_1$ with supremum
$\gamma$, starting with $\gamma_0^*= \gamma_0$. 
By induction choose $r_i\in \App\restrictedto \gamma_i^*$ for
$i<\omega_1$ such that: 
$$r_i \leend r_j \hbox{ for $i < j<\omega_1$};\leqno(1)$$
$$q_j\restrictedto \gamma_i^*\le r_i \hbox{ for $j<\delta$ and $i<
\omega_1$}.\leqno (2)$$
Here for each $i$ the inductive hypothesis is applied to $q_j\restrictedto
\gamma_i^*$, $r_i$, and $\gamma_i$. 

\noindent {\sl The case $\gamma$ a limit of cofinality $\aleph_0$.} 


End of Garbage Heap

\vfill\eject		


\appendix

\proclaim{Omitting types} \endproclaim

In \S1 we made (implicit) use of the combinatorial
principle developed in [ShHL162].  
In the context of this paper, 
this is a combinatorial refinement of forcing with
$\App$, which gives (in the ground model)
a $\Pp_3$-name $\name \FF$ for a filter with the required properties
in a $\Pp_3$-generic extension.  We now review this material.  Our
discussion
overlaps with the discussion in [Sh326], but will be more complete in
some technical 
respects and less complete in others.  
We begin in
sections A1-A5 
by presenting the material of [Sh162] as it was summarized in [Sh326].
However the setup of [Sh162] can be (and should be) 
tailored more closely to the
applications, and we will present a second setup which is more
convenient in sections A6-A10. One could take the view that the
axioms given in section A6 below should supercede the axioms given in
section A1, and one should check that the proofs of [Sh162] work with
these new axioms.  Since this would be awkward in practice, we  take a
different route, showing that the two formalisms are equivalent.  

After dealing with this technical point, we will
not explain in any more detail the way this principle is applied, as
that aspect is dealt with at 
great length in a very similar context in [Sh326].  For the reader who
is not familiar with [Sh162] the discussion in the appendix to [Sh326]
should be more useful than the present discussion.

\\ Uniform partial orders

We review the formalism of [Sh162].

With the cardinal $\lambda$ fixed, a partially ordered set $(\PP,<)$ is
said to be {\sl standard $\lambda^+$-uniform} if $\PP\includedin
\lambda^+\times P_{\lambda}(\lambda^+)$ (we refer here to subsets of
$\lambda^+$ of size strictly less than $\lambda$), has the
following properties (if  $p=(\alpha,u)$ we write
$\dom p$ for $u$, and we write $\PP_\alpha$ for $\{p\in \PP:\dom p
\includedin \alpha\}$):

\item {1.} If $p\le q$ then $\dom p\includedin \dom q$.

\item {2.} For all $p \in \PP$ and $\alpha < \lambda^+$ there exists a
$q \in \PP$ with $q \le p$ and $\dom q = \dom p \intersect \alpha$;
furthermore, there is a unique maximal such $q$, for which we write
$q=p \restrictedto\alpha$.

\item {3.}
{\sl $($Indiscernibility$)$} If $p=(\alpha, v) \in \PP$ and $h:v
\rightarrow v'\includedin \lambda^+$ is an order-isomorphism onto $V'$ then
$(\alpha, v') \in
\PP$.  We write $h[p]=(\alpha, h[v])$.
Moreover, if $q\le p$ then $h[q]\le
h[p]$.

\item {4.}
{\sl $($Amalgamation$)$} For every $p,q \in \PP$ and $\alpha<\lambda^+$,
if $p\restrictedto\alpha \le q$ and $\dom p \intersect \dom q= \dom p
\intersect \alpha$, then there exists $r \in \PP$ so that $p,q \le r$.

\item {5.}
For all $p,q,r \in \PP$ with $p,q \le r$ there is $r' \in \PP$ so that
$p,q \le r'$ and $\dom r'=\dom p \union
\dom q$.

\item {6.}
If $(p_i)_{i < \delta} $ is an increasing sequence of length less than
$\lambda$, then it has a least upper bound $q$, with domain \
$\Union_{i<\delta} \dom p_i$; we will write $q=\Union_{i<\delta}
p_i$, or more succinctly: $q=p_{<\delta}$.

\item {7.}
For limit ordinals $\delta$, $p \restrictedto \delta =\Union_{\alpha <
\delta} p \restrictedto\alpha$.

\item {8.}
If $(p_i)_{i < \delta}$ is an increasing sequence of length less than
$\lambda$, then  $(\Union_{i < \delta} p_i)\restrictedto\alpha\ =
\  
\Union_{i < \delta}
(p_i \restrictedto \alpha).$

It is shown in [ShHL162] that under a diamond-like
hypothesis, such partial orders admit reasonably generic objects.  The
precise formulation is given in A5 below.

\\ Density systems

Let $\PP$ be a standard $\lambda^+$-uniform partial order.  For
$\alpha<\lambda^+$, $\PP_\alpha$ denotes the restriction of $\PP$ to
$p\in \PP$ with domain contained in $\alpha$.  A subset $ G$ of
$\PP_\alpha$ is an {\sl admissible ideal} (of $\PP_\alpha$)
 if it is closed downward, is
$\lambda$-directed (i.e.  has upper bounds for all small subsets), and
has no proper directed extension within $\PP_\alpha$.  For $G$ an
admissible ideal in $\PP_\alpha$, $\PP/G$ denotes the restriction of
$\PP$ to $\{p\in \PP: p\restrictedto \alpha\in G\}$.

If $G$ is an admissible ideal in $\PP_\alpha$ and
$\alpha<\beta<\lambda^+$, then an $(\alpha,\beta)$-{\sl density system}
for $G$ is a function $D$ from pairs $(u,v)$ in $P_\lambda(\lambda^+)$
with $u\includedin v$ into subsets of $\PP$ with the following
properties: \item{(i)} $D(u,v)$ is an upward-closed dense subset of
$\{p\in \PP/G:\dom p\includedin v\union \beta\}$; 
\item{(ii)} For pairs
$(u_1,v_1),(u_2,v_2)$ in the domain of $D$, if
$u_1\intersect\beta=u_2\intersect\beta$ and
$v_1\intersect\beta=v_2\intersect\beta$, and there is an order
isomorphism from $v_1$ to $v_2$ carrying $u_1$ to $u_2$, then for any
$\gamma$ we have $(\gamma,v_1)\in D(u_1,v_1)$ iff $(\gamma,v_2)\in
D(u_2,v_2)$.
\medskip

An admissible ideal $G'$ (of $\PP_\gamma$) is said to
{\sl meet} the $(\alpha,{\beta})$-density system $D$ for $G$ if
$\gamma\ge\alpha$,
$G'\ge G$
and for each $u\in P_\lambda(\gamma)$ there is $v\in P_\lambda(\gamma)$
containing $u$ such that $G'$ meets $D(u,v)$.

\\ The genericity game

Given a standard $\lambda^+$-uniform partial order $\PP$, the {\sl
genericity game} for $\PP$ is a game of length $\lambda^+$ played by
Guelfs and Ghibellines, with Guelfs moving first.  The Ghibellines
build an increasing sequence of admissible ideals meeting density
systems set by the Guelfs.  Consider stage $\alpha$.  If $\alpha$ is a
successor, we write $\alpha^-$ for the predecessor of $\alpha$;  if
$\alpha$ is a limit, we let $\alpha^-=\alpha$. Now at stage $\alpha$  
for every $\beta<\alpha$ an admissible ideal $G_\beta$
in some $\PP_{\beta'}$ is given, and one can
check that there is a unique admissible ideal $G_{\alpha^-}$ in
$\PP_{\alpha^-}$ containing
$\Union_{\beta<\alpha}G_{\beta'}$ (remember A 3.1(5)) 
[Lemma 1.3, ShHL 162].  The Guelfs now
supply at most $\lambda$\
density systems $D_i$  over $G_{\alpha^-}$ for
$(\alpha,\beta_i)$ and also fix an element $g_\alpha$
in $\PP/G_\alpha^-$.
Let $\alpha'$ be minimal such that $g_\alpha\in \PP_{\alpha'}$ and
$\alpha'\ge\sup\,\beta_i$.  The Ghibellines then build an admissible ideal
$G_{\alpha'}$ for $\PP_{\alpha'}$ containing $G_\alpha^-$ as well as
$g_\alpha$, and meeting all specified density systems, or forfeit the
match; they let $G_{\alpha''}=G_{\alpha'}\cap \alpha''$ when
 $\alpha\le \alpha''<\alpha'$.
   The main result is that the Ghibellines can win with a little
combinatorial help in predicting their opponents' plans.

For notational simplicity, we assume that $G_\delta$ is an
$\aleph_2$-generic ideal on $\App\restrictedto \delta$,
when $\cf \delta=\aleph_2 $, which is true
on a club in any case.

\\  $Dl_\lambda$

The combinatorial principle $\rm Dl_\lambda$ states that there are subsets
$Q_\alpha$ of the power set of $\alpha$ for $\alpha<\lambda$ such that
$|Q_\alpha|<\lambda$, and for any $A\includedin \lambda$ the set
$\{\alpha:A\intersect\alpha\in Q_\alpha\}$ is stationary.  This follows
from $\diamond_\lambda$ or inaccessibility, obviously, and Kunen showed
that for successors, $\rm Dl$ and $\diamond$ are equivalent.  In addition
$Dl_\lambda$ implies $\lambda^{<\lambda}=\lambda$.  

\bigskip

\\ A general principle

\proclaim{Theorem} 

Assuming $Dl_\lambda$, the Ghibellines can win any standard
$\lambda^+$-uniform $\PP$-game. 
\endproclaim

This is Theorem 1.9 of [ShHL 162].


\\ Uniform partial orders revisited

We introduce a second formalism that fits the setups encountered in
practice more closely.  In our second version we write
``quasiuniform'' rather than ``uniform'' throughout as the axioms have
been weakened slightly.

With the cardinal $\lambda$ fixed, a partially ordered set $(\PP,<)$ is
said to be {\sl standard $\lambda^+$-quasiuniform} if $\PP\includedin
\lambda^+\times P_{\lambda}(\lambda^+)$  has the
following properties 
(if  $p=(\alpha,u)$ we write
$\dom p$ for $u$, and we write $\PP_\alpha$ for $\{p\in \PP:\dom p
\includedin \alpha\}$):

\item {$1'$.} If $p\le q$ then $\dom p\includedin \dom q$.

\item {$2'$.} For all $p \in \PP$ and $\alpha < \lambda^+$ there exists a
$q \in \PP$ with $q \le p$ and $\dom q = \dom p \intersect \alpha$;
furthermore, there is a unique maximal such $q$, for which we write
$q=p \restrictedto\alpha$.

\item {$3'$.}
{\sl (Indiscernibility)} If $p=(\alpha, v) \in \PP$ and $h:v
\rightarrow v'\includedin \lambda^+$ is an order-isomorphism onto $V'$ then
$(\alpha, v') \in
\PP$.  We write $h[p]=(\alpha, h[v])$.
Moreover, if $q\le p$ then $h[q]\le
h[p]$.

\item {$4'$.}
{\sl (Amalgamation)} For every $p,q \in \PP$ and $\alpha<\lambda^+$,
if $p\restrictedto\alpha \le q$ and $\dom p \intersect \dom q= \dom p
\intersect \alpha$, then there exists $r \in \PP$ so that $p,q \le r$.

\item {$5'$.}
If $(p_i)_{i < \delta} $ is an increasing sequence of length less than
$\lambda$, then it has an upper bound $q$.

\item {$6'$.}
If $(p_i)_{i < \delta}$ is an increasing sequence of length less than
$\lambda$ of members of $\PP_{\beta+1}$, with $\beta <\lambda^+$ and
if $q\in \PP_\beta$ satisfies $p_i\restrictedto \beta\le q$ for all
$i< \delta$, then $\{p_i:i < \delta\}\union \{q\}$ has an upper bound
in $\PP$. 

\item {$7'$.}
If $(\beta_i)_{i < \delta}$ is a strictly increasing sequence of
length less than
$\lambda$, with each $\beta_i<\lambda^+$, and $q\in \PP$,
$p_i\in \PP_{\beta_i}$, 
with $q\restrictedto \beta_i \le p_i$, 
then $\{p_i:i< \delta\}\union \{q\}$ has an
upper bound.

\item {$8'$.}
Suppose $\xi, \zeta$ are limit ordinals less than $\lambda$, and
$(\beta_i)_{i< \zeta}$ is a strictly increasing continuous sequence of
ordinals less than $\lambda^+$. Let 
$I(\xi,\zeta) := (\zeta+1)\times (\xi +1) -\{(\zeta,\xi)\}$. 
Suppose that for $(i,j)\in I(\xi,\zeta)$
we have $p_{ij}\in
\PP_{\beta_i}$ such that 
$$
\eqalign{
i\le i' &\implies p_{ij}=p_{i'j}\restrictedto \beta_i;\cr
j\le j' &\implies p_{ij}\le p_{ij'}\cr}
$$

Then $\{p_{ij}: (i,j)\in I(\xi,\zeta)\}$ has an upper bound in
$\PP$.  

\\ Density systems revisited

Let $\PP$ be a standard $\lambda^+$-quasiuniform partial order.  
A subset $ G$ of
$\PP_\alpha$ is a {\sl quasiadmissible ideal} (of $\PP_\alpha$)
if it is closed downward and is 
$\lambda$-directed (i.e.  has upper bounds for all small subsets).
For $G$ a quasiadmissible ideal in $\PP_\alpha$, $\PP/G$ denotes the
restriction of 
$\PP$ to $\{p\in \PP: p\restrictedto \alpha\in G\}$.

If $G$ is a quasi-admissible ideal in $\PP_\alpha$ and
$\alpha<\beta<\lambda^+$, then an $(\alpha,\beta)$-{\sl density
system} for $G$ is a function $D$ from sets $u$ in
$P_\lambda(\lambda^+)$ 
into subsets of $\PP$ with the following
properties: 
\item{(i)} $D(u)$ is an upward-closed dense subset of
$\PP/G$;
\item{(ii)} For pairs
$(u_1,v_1)$ and $(u_2,v_2)$ with $u_1$, $u_2$ 
in the domain of $D$, and $v_1,v_2\in
P_\lambda(\lambda^+)$ with $u_1\includedin v_1$, $u_2\includedin v_2$,
if  
$u_1\intersect\beta=u_2\intersect\beta$ and
$v_1\intersect\beta=v_2\intersect\beta$, and there is an order
isomorphism from $v_1$ to $v_2$ carrying $u_1$ to $u_2$, then for any
$\gamma$ we have $(\gamma,v_1)\in D(u_1)$ iff $(\gamma,v_2)\in
D(u_2)$.
\medskip

For $\gamma \ge \alpha$, 
a quasiadmissible ideal $G'$ of $\PP_\gamma$ is said to
{\sl meet} the $(\alpha,{\beta})$-density system $D$ for $G$ if
$G'\ge G$
and for each $u\in P_\lambda(\gamma)$ $G'$ meets $D(u,v)$.

\\ The genericity game revisited

Given a standard $\lambda^+$-quasiuniform partial order $\PP$, the {\sl
genericity game} for $\PP$ is a game of length $\lambda^+$ played by
Guelfs and Ghibellines, with Guelfs moving first.  The Ghibellines build
an increasing sequence of admissible ideals meeting density systems set
by the Guelfs.  Consider stage $\alpha$.  If $\alpha$ is a successor,
we write $\alpha^-$ for the predecessor of $\alpha$;  if $\alpha$ is a
limit, we let $\alpha^-=\alpha$. Now at stage $\alpha$ 
 for every $\beta<\alpha$ an admissible ideal $G_\beta$
 in some $\PP_{\beta'}$ is given.  The Guelfs now
supply at most $\lambda$\
density systems $D_i$  over $G_{\alpha^-}$ for
$(\alpha,\beta_i)$ and also fix an element $g_\alpha$
in $\PP/G_\alpha^-$.
Let $\alpha'$ be minimal such that $g_\alpha\in \PP_{\alpha'}$ and
$\alpha'\ge\sup\,\beta_i$.  The Ghibellines then build an admissible ideal
$G_{\alpha'}$ for $\PP_{\alpha'}$ containing $\Union_{\beta < \alpha}
G_\beta$ as well as 
$g_\alpha$, and meeting all specified density systems, or forfeit the
match; they let $G_{\alpha''}=G_{\alpha'}\cap \alpha''$ when
 $\alpha\le \alpha''<\alpha'$.
   The main result is that the Ghibellines can win with a little
combinatorial help in predicting their opponents' plans.

\proclaim{Theorem} 

Assuming $Dl_\lambda$, the Ghibellines can win any standard
$\lambda^+$-uniform $\PP$-game. 
\endproclaim

We will show this is equivalent to the version given in [ShHL162]. 

\\ The translation

To match up the uniform and quasiuniform settings, we give a
translation of the quasiuniform setting back into the uniform setting;
there is then an accompanying translation of density systems and of
the genericity game. 
So we assume that the standard $\lambda^+$-quasiuniform partial order
$\PP$ is given and we will define an associated partial ordering
$\PP'$. 

The set of elements of $\PP'$ is the set of sequences $\pbar =
(p_{ij}, \beta_i)_{i<\zeta,j<\xi}$ such that:\
$$\def\offsetsize{37.06 true pt}
\eqalign{
\offset{(a)}&\hbox{$\zeta,\xi < \lambda$; $\beta_i$ is
strictly increasing;}\cr
\offset{(b)}&\hbox{$p_{ij}=p_{i'j}\restrictedto \beta_i$, and $\beta_i\in \dom p_{i'j}$,  for $i <
i'$;}\cr 
\offset{(c)}&\hbox{$p_{ij}< p_{ij'}$ for $j < j'$;}\cr
\offset{(d)}&\hbox{If $\alpha = \delta+\alpha'\in \dom p_{ij}$ with
$\alpha'< \lambda$ and $\delta$ is divisible by $\lambda$ and of
cofinality less than $\lambda$, then \qquad
}\cr 
&\hbox{$\delta\intersect \dom p_{ij}$ 
is unbounded in $
\delta$.}\cr}
$$

For $\pbar\in \PP'$ let 
$\dom \pbar = 
\{ \delta+n: \exists i,j\, \dom
p_{ij}\intersect [(\delta+\epsilon_\delta+n)\lambda,
(\delta+\epsilon_\delta+n+1)\lambda) \ne \emptyset\}$, 
where $\delta$ is a limit ordinal or $0$ and where $\epsilon_\delta$
is $0$ if $\cof \delta$ is $\lambda$, and is $1$ otherwise.
We can represent the elements of $\PP'$ naturally by codes of the type
used in \S A1, so that the domain as defined here is the domain in the
sense of this coding as well.

Now we define the order on $\PP'$. For $\pbar,\qbar\in \PP'$ we have
the associated ordinals (such as $\zeta^\qbar$), and the elements
$p_{ij}, q_{ij}$ of $\PP$.  We say $\pbar \le \qbar$
if one of the following occurs:

\item{1.} $\pbar=\qbar$;
\item{2.} $\zeta^\pbar=\zeta^\qbar$,
$\beta_i^\pbar=\beta_i^\qbar$ 
for $i< \zeta^\pbar$, and there is $j'< \xi^\qbar$
such that $p_{ij}\le q_{ij'}$ for all $i< \zeta^\pbar$ and
$j< \xi^\pbar$.  
\item{3.}$\xi^\pbar=\xi^\qbar$ and 
there is $i'<\zeta^\qbar$ such that $p_{ij}\le q_{i_jj}$ for all
$i<\zeta^\pbar$ and $j<\xi^\pbar$. 
\item{4.} There are $i'$, $j'$ such that $p_{ij}\le q_{i'j'}$ for all
$i< \zeta^\pbar$ and $j< \xi^\qbar$.

The first thing to be checked is that this is transitive.
We will refer to relations of the type described in (2-4) above as
{\sl vertical}, {\sl horizontal}, or {\sl planar} respectively.
The equality relation may be considered as being of all three types.
With regard to transitivity, 
if $\pbar\le \qbar\le \rbar$,  
then if both of the inequalities involved are
horizontal, or both are 
vertical, we have an inequality $\pbar \le \rbar$
of the same type; and otherwise we have a planar inequality $\pbar \le
\rbar$. 

We do not insist on asymmetry; if one wishes to have a partial order
in the strict sense then it will be necessary to factor out an
equivalence relation.

\\ Properties (A1.1-4)

We claim that if $\PP$ is a partial order with properties 
$1'$-$8'$ of \S A6, then the associated partial ordering $\PP'$ enjoys
properties $1$-$8$ of \S A1. 
The first four properties were assumed for $\PP$; we have to check
that they are retained by $\PP'$. 

\item {1.} 
If $\pbar\le \qbar$ then $\dom \pbar\includedin \dom \qbar$. 

\proof:
If $\pbar \le \qbar$ then $\Union \dom p_{ij}\le \Union \dom q_{i'j'}$
by (1) applied to $\PP$ and hence $(1)$ holds for $\PP'$ by applying
the definition of $\dom$ in $\PP'$. 
\qed

\item {2.} 
For all $\pbar \in \PP'$ and $\alpha < \lambda^+$ there exists a 
$\qbar \in \PP'$ with 
$\qbar \le \pbar$ and $\dom \qbar = \dom \pbar \intersect \alpha$; 
furthermore, there is a unique maximal such $\qbar$, for which we
write $\qbar=\pbar \restrictedto\alpha$.

\proof:
Let $\alpha'=\alpha\cdot \lambda$, $\zeta'=\{i:\beta_i^\pbar <
\alpha'\}$, and
$p_{ij}'=p_{ij}\restrictedto
\alpha'$ for $i<\zeta'$. Set $\pbar\restrictedto \alpha =
(p_{ij}',\beta_i)_{i< \zeta', j< \xi ^\pbar}$. 
\vfill \eject

\item {3.}
{\sl $($Indiscernibility$)$} 
If $\pbar=(\alpha, v) \in \PP'$ and 
$h:v \rightarrow v'\includedin \lambda^+$ is an order-isomorphism onto
$V'$ then $(\alpha, v') \in \PP'$.  
We write $h[\pbar]=(\alpha, h[v])$.
Moreover, if $\qbar\le \pbar$ then $h[\qbar]\le h[\pbar]$.

\bigskip

\item {4.}
{\sl $($Amalgamation$)$} 
For every $\pbar,\qbar \in \PP'$ and $\alpha<\lambda^+$, if
$\pbar\restrictedto\alpha \le \qbar$ and 
$\dom \pbar \intersect \dom \qbar= \dom \pbar \intersect \alpha$, then
there exists $r \in \PP'$ so that $\pbar,\qbar \le r$. 

\bigskip

\\ Property (A1.5)

We consider the fifth property:

\item {5.}
For all $\pbar,\qbar,\rbar \in \PP'$ with $\pbar,\qbar \le r$ there is
$r' \in \PP'$ so that $\pbar,\qbar \le r'$ and 
$\dom r'=\dom \pbar \union \dom \qbar$.

\bigskip

\\ Properties (A1.6-8)

The last three properties are:

\item {6.}
If $(\pbar_i)_{i < \delta} $ is an increasing sequence of length less
than $\lambda$, then it has a least upper bound $\qbar$, with domain 
$\Union_{i<\delta} \dom \pbar_i$;  
we will write $\qbar=\Union_{i<\delta} \pbar_i$, 
or more succinctly: $\qbar=\pbar_{<\delta}$.

\bigskip

\item {7.}
For limit ordinals $\delta$, 
$\pbar \restrictedto \delta = \Union_{\alpha < \delta} 
\pbar \restrictedto\alpha$. 

\bigskip

\item {8.}
If $(\pbar_i)_{i < \delta}$ is an increasing sequence of length less
than $\lambda$, then 
$(\Union_{i < \delta} \pbar_i)\restrictedto\alpha = 
\Union_{i < \delta} (\pbar_i \restrictedto \alpha).$

\bigskip

\\ Application

In our application we identify $\App$ with a standard
$\aleph_2^+$-uniform partial order via a certain coding. We first
indicate a natural coding which is not quite the right one, then
repair it.
\medskip

\noindent {\bf First Try}

An approximation $q=(\AA,\name \FF,\name{\fakebf\epsilon})$ will be
identified with a pair $(\tau,u)$, where $u=\AA$, and $\tau$ is the
image of $q$ under the canonical order-preserving map $h:\AA\with
\otp(\AA)$.  
One important point is that the first parameter $\tau$ comes from a
fixed set $T$ of size $2^{\aleph_1}=\aleph_2$; so if we enumerate $T$
as $(\tau_\alpha)_{\alpha<\aleph_2}$ then we can code the pair
$(\tau_\alpha,u)$ by the pair $(\alpha,u)$.  Under these successive
identifications, $\App$ becomes a standard $\aleph_2^+$-uniform
partial order, as defined in \S A1.  Properties1 , 2, 4, 5, and 6 are
clear, as is 7, in view of the uniformity in the iterated forcing
$\Pp$, and properties 3, 8 were, stated in 1.7 and 1.9 (4).
\comment{This part will change}

The difficulty with this approach is that in this formalism, density
systems cannot express nontrivial information: any generic ideal meets
any density system, because for $q\le q'$ with $\dom q=\dom q'$, we
will have $q=q'$; thus $D(u,u)$ will consist of all $q$ with $\dom q =
u$, for any density system $D$.

So to recode $\App$ in a way that allows nontrivial density systems to
be defined, we proceed as follows.
\medskip

\noindent{\bf Second Try}

Let $\iota:\aleph_2^+\with \aleph_2^+\times \aleph_2$ be order
preserving where $\aleph_2^+ \times \aleph_2$ is ordered
lexicographically.  Let $\pi:\aleph_2^+\times \aleph_2\to \aleph_2^+$
be the projection on the first coordinate. First encode $q$ by
$\iota[q]=(\iota[\AA],\ldots)$, then encode $\iota[q]$ by
$(\tau,\pi[\AA])$, where $\tau$ is defined much as in the first try --
a description of the result of collapsing $q$ into
$\otp\pi[\AA]\times\aleph_2$, after which $\tau$ is encoded by an
ordinal label below $\aleph_2$. The point of this is that now the
domain of $q$ is the set $\pi[\AA]$, and $q$ has many extensions with
the same domain.  After this recoding, $\App$ again becomes a
$\aleph_2^+$-uniform partial ordering, as before.  We will need some
additional notation in connection with the indiscernibility condition.
It will be convenient to view $\App$ simultaneously from an encoded
and a decoded point of view.  One should now think of $q\in \App$ as a
quadruple $(u,\AA,\name\FF, {\fakebf\epsilon})$ with $\AA\includedin
u\times \aleph_2$.  If $h:u\with v$ is an order isomorphism, and $q$
is an approximation with domain $u$, we extend $h$ to a function $h_*$
defined on $\AA^q$ by letting it act as the identity on the second
coordinate.  Then $h[q]$ is the transform of $q$ using $h_*$, and has
domain $v$.

\bigskip\bigskip
For notational simplicity, we assume that $G_\delta$ is an
$\aleph_2$-generic ideal on $\App\restrictedto \delta$, when $\cf
\delta=\aleph_2 $ which is true on a club in any case.
\comment{Does this remark go anywhere?}

\vfill\eject		

 
\references
\setrefwidth{ShHL 162}

\ref AxKo. J. Ax and S. Kochen, {\sl Diophantine problems over local
rings I.}, Amer. J. Math. {\bf 87} (1965), 605--630.

\ref Ch. G. Cherlin, {\sl Ideals of integers in nonstandard number
fields}, in: {\bf Model Theory and Algebra}, LNM 498, Springer, New 
York, 1975, 60--90.

\ref DW. H. G. Dales and W. H. Woodin, {\bf An Introduction to 
Independence for Analysts}, CUP Cambridge, 1982.

\ref Ke. H. J. Keisler, {\sl Ultraproducts which are not saturated},
J. Symbolic Logic {\bf 32} (1967), 23--46.

\ref KeSc. H. J. Keisler and J. Schmerl, {\sl }, \JSL 56 (1991), {\sl
Making the hyperreal line both saturated and complete}, 1016--1025.


\ref LLS. Ronnie Levy, Philippe Loustaunau, and Jay Shapiro, {\sl The
prime spectrum of an infinite product of copies of $\Zz$}, Fund.
Mathematicae {\bf 138} (1991), 155--164.   

\ref Mo. J. Moloney,  {\sl Residue class domains of the ring of
convergent sequences and of $C^{\infty}([0,1],\Rr)$}, Pacific J. Math.
{\bf 143} (1990), 1--73. 

\ref Ri. D. Richard, {\sl De la structure additive \`a la saturation des
mod\'eles de Peano et \`a une classification des sous-langages de
l'arithm\'etique}, in {\bf Model Theory and Arithmetic (Paris,
1979/80)}, C. Berline et al eds., LNM 890, Springer, New York, 1981,
270--296.

\ref Sh-a. S. Shelah, {\bf Classification Theory and the Number of
Non-Isomorphic 
Models}, {\sl North Holland Publ. Co.},
Studies in Logic and the foundation of Math., vol. 92, 1978.

\ref Sh-c.  \same, {\bf Classification Theory and the Number of
Non-isomorphic Models}, 
revised, {\sl North Holland Publ. Co.}, Studies in Logic 
and the foundation of Math., Vol.92, 1990, 705+xxxiv.

\ref Sh72. \same, {\sl Models with
second order properties I. Boolean algebras with no undefinable
automorphisms}, Annals Math. Logic {\bf 14} (1978), 57--72.

\ref Sh107. \same, {\sl Models with
second order properties IV. A general method and eliminating
diamonds}, Annals Math. Logic {\bf 25} (1983), 183--212.

\ref ShHL162. \same, B. Hart, and C. Laflamme, {\sl Models with
second order properties V. A general principle}, Annals Pure Applied
Logic, to appear. 

\ref Sh326. \same, {\sl Vive la diff\'erence I. Nonisomorphism of
ultrapowers of countable models}, Proceedings of the Oct. 1989 MSRI
Conference on Set Theory, J. Judah, W. Just, and W. H. Woodin eds., to
appear. 

\ref Sh482. In preparation.
\end
\vfill\eject		
\nopagenumbers
\end

\bye